\documentclass[11pt]{article}

\usepackage{amsmath}
\usepackage{dsfont}
\usepackage{amsfonts}
\usepackage{amssymb}
\usepackage{mathtools}
\usepackage{theorem}
\usepackage{ifthen}
\usepackage{mathrsfs}
\usepackage{cancel}
\usepackage{soul}
\usepackage{cite}
\usepackage{color}
\usepackage{hyperref}
\definecolor{links}{rgb}{0.36,0.54,0.66}
\hypersetup{
   colorlinks = true,
    linkcolor = black,
     urlcolor = blue,
    citecolor = blue,
    filecolor = blue,
   }
\usepackage{tikz}
\usepackage{comment}

\newcommand{\eps}{\varepsilon}
\newcommand {\R}{\mathbb R}
\newcommand{\qed}{\hfill\raisebox{1.2mm}{\fbox{}}}
\newcommand{\omitted}{\hfill\triangle}
\newcommand{\bdry}{\partial}
\newcommand{\card}{\#}

\newboolean{comments}

\setboolean{comments}{true}

\newcommand{\dl}[1]{\ifthenelse{\boolean{comments}}{\textcolor{blue}{\small #1}}{}}
\newcommand{\dlnew}[1]{\ifthenelse{\boolean{comments}}{\textcolor{teal}{\small #1}}{}}
\newcommand{\ak}[1]{\ifthenelse{\boolean{comments}}{\textcolor{red}{\small #1}}{}}
\newcommand{\aknew}[1]{\ifthenelse{\boolean{comments}}{\textcolor{magenta}{\small #1}}{}}

\newtheorem{Theorem}{Theorem}[section]
{\theorembodyfont{\rmfamily} 
\newtheorem{Lemma}[Theorem]{Lemma}
\newtheorem{Proposition}[Theorem]{Proposition}
\newtheorem{Corollary}[Theorem]{Corollary}
{\theorembodyfont{\rmfamily}
\newtheorem{Remark}[Theorem]{Remark}}
{\theorembodyfont{\rmfamily} \newtheorem{Definition}{Definition}}
{\theorembodyfont{\rmfamily} \newtheorem{Example}{Example}}

\parskip=5pt plus 1pt minus 1pt

\title{\bf Robust signal decompositions on the circle}

\author{Aral Köse\thanks{A. Köse is with Bo\u{g}azi\c{c}i University, Turkey (email: aral.kose@std.bogazici.edu.tr).}
\ and Daniel Liberzon\thanks{D. Liberzon is with the Coordinated Science Laboratory, University of
	Illinois Urbana-Champaign, Urbana, IL 61801 USA (email: liberzon@uiuc.edu). His work was supported in part by the NSF CMMI-2106043 and AFOSR MURI FA9550-23-1-0337 grants.
 }}

\begin{document}

\maketitle

\begin{abstract} We consider the problem of decomposing a piecewise constant function on the circle into a sum of indicator functions of closed circular disks in the plane, whose number and location are not a priori known. This represents a situation where an agent moving on the circle is able to sense its proximity to some landmarks, and the goal is to estimate the number of these landmarks and their possible locations---which can in turn enable control tasks such as motion planning and obstacle avoidance. Moreover, the exact values of the function at its discontinuities (which correspond to disk boundaries for the individual indicator functions) are not assumed to be known to the agent. We introduce suitable notions of robustness and degrees of freedom
to single out those decompositions that are more desirable, or more likely, given this non-precise data collected by the agent. We provide a characterization of robust decompositions and give a procedure for generating all such decompositions. When the given function admits a robust decomposition, we compute the number of possible robust decompositions and derive bounds for the number of decompositions maximizing the degrees of freedom.
\end{abstract}

\section{Introduction}

Imagine an agent moving along a circular path in the plane with some stationary landmarks, whose number and exact locations are unknown to the agent. Suppose that each landmark transmits an omnidirectional signal with a finite range, which we can model as a function that equals~1 inside a circular disk centered at the landmark and~0 outside. The boundaries of these disks, whose radii are in general different, may intersect the agent's path at one or two points or not at all. As the agent moves along its path, it can perceive these signals and so it knows, at each point, the number of landmarks that are within range. It cannot, however, identify different landmarks by their signals, and neither can it discern anything about each signal's strength other than its presence or absence. The agent's knowledge of its position on the circle may also not be precise, and the signal transmissions or measurements may occur with some sampling frequency rather than continuously in time. For these reasons, all that the agent can reliably reconstruct is a sequence of nonnegative integers corresponding to local landmark counts around the circle, and it may not be sure of the precise count at the exact points where this count changes.

In this scenario, we want to pose the following questions: Can the agent figure out the total number of landmarks (excluding, of course, those whose signals do not reach any points on the circle)? Can it reconstruct some qualitative information about how these landmarks---or, more precisely, the disks around them where their signals equal 1---are positioned relative to the circle; i.e., do they intersect it at a single point or along an arc or cover it entirely? Moreover, in view of what the agent is able to measure, we may naturally prefer some landmark configurations to others. Namely, we will single out those landmark arrangements which, even when slightly perturbed, would still give the same sequence of the agent's local landmark counts around the circle. We will term such landmark arrangements \emph{robust} and argue that by being ``in general position" they provide more likely explanations of the data collected by the agent, compared to ``special" ones for which the agent's data would change under arbitrarily small perturbations to the landmark positions or the ranges of their signals. We can then ask, how many such robust landmark configurations are possible, and what do they qualitatively look like? Interestingly, we will see that when robustness is imposed, the agent's lack of knowledge of landmark counts at some isolated points is in fact inconsequential---these missing values can be reconstructed.

The setting described above shares some aspects with the well-known problem of simultaneous localization and mapping (SLAM) in robotics; see, e.g., \cite{dissanayake2001solution,thrun2002probabilistic,lavalle-book} and the references therein. In SLAM, the agent's task is to reconstruct a map of an uncertain environment---typically containing some landmarks---and to localize itself within this environment. Our emphasis here is on the mapping rather than on the localization. Among the various versions of SLAM appearing in the literature, the one considered in the recent work~\cite{LM-SLAM-hscc2025} has some common features with our set-up and includes, in particular, an example where an agent moving on a circle is able to reconstruct its position based on signal measurements (of the same on/off type as above) from a single landmark. Multiple landmarks are not addressed in that work, and its general emphasis is on the localization (more precisely, on describing its limitations by characterizing agent's states that cannot be distinguished from one another).

The paper~\cite{euler} studies a problem very closely related to ours, although motivated from a slightly different point of view. Instead of a moving agent, the authors of~\cite{euler} envision a network of sensors, each able to count the number of landmarks in its vicinity. They then ask how the sensors can merge their local counts into a global one. They assume that precise local counts are known at all points and, for topological reasons, they do not allow the presence of landmarks whose signals reach all sensors. Nevertheless, as we explain in more detail in Section~\ref{s-generating-counting} below, it is possible to reduce our case to theirs by appropriate pre-processing (essentially, by removing such landmarks before counting) which allows us to recover some of our results from theirs. Also, while we restrict the agent's path and the supports of the landmarks' signals to be circular, in~\cite{euler} both can have more general shapes. On the other hand, \cite{euler} does not consider the robustness property, which is the main focus of the present paper.

Some brief highlights of our contributions are as follows. In Section~\ref{s-prelim} we formalize the notion of a \emph{decomposition} of a piecewise constant function representing the agent's local landmark counts into a sum of signal functions from individual landmarks (Definition~\ref{d-decomposition}) and introduce equivalence classes of sequences of these local landmark counts with respect to rotations around the circle. In Section~\ref{s-robust} we define \emph{robust} decompositions to be those whose sequence representations are invariant under sufficiently small perturbations (Definition~\ref{d-robust}). We then establish a necessary and sufficient condition for robustness, which says that when the boundary of the support disk of a landmark's signal intersects the circle it should do so at two distinct points, and these points of intersection should be different for different landmarks (Proposition~\ref{p-robust}). In Section~\ref{s-top-prop} we define a larger class of decompositions, which we call \emph{properly upper semicontinuous} (properly u.s.c.) and which include all robust decompositions. We provide an equivalent characterization of this property which extends the corresponding result for robust decompositions (Lemma~\ref{l-pusc-char}) and then derive a lower bound on how many signals within a properly u.s.c.\ decomposition must be identically~1 on the circle (Proposition~\ref{p-baseline-levels}) and an exact formula for how many other signals, covering a smaller arc of the circle, there must be (Corollary~\ref{c-intervals}). In Section~\ref{s-dof} we examine properly u.s.c.\ and robust decompositions from a different angle by introducing \emph{degrees of freedom}, which tell us in how many dimensions we can perturb the landmark positions and their signal ranges to end up with exactly the same total signal function. We show that this number of degrees of freedom is maximized by properly u.s.c.\ decompositions satisfying an additional condition, which includes all robust decompositions (Proposition~\ref{p-pusc-max-dof}). For such degrees-of-freedom-maximizing decompositions, in Section~\ref{s-generating-counting} we derive lower and upper bounds on the total number of landmarks that they can contain (Proposition~\ref{p-number-of-signals-pusc}) and determine how many distinct decompositions of this kind exist (Proposition~\ref{p-counting}). Detailed examples are included to illustrate the results.

\section{Signal functions and decompositions}\label{s-prelim}

Suppose we start with the data $(m_i,r_i)$, $i=1,\dots, n$, where for each $i$, $m_i\in\R^2$ is the location of the $i$th signal source (landmark) and $r_i$ is the signal's radius. The individual \textit{signal functions} are $h_i:\R^2\to\{0,1\}$ given by
\setlength{\abovedisplayskip}{3pt}\setlength{\belowdisplayskip}{3pt}
\begin{equation}\label{e-h-binary-def}
h_i(x):=\left\{\begin{array}{ll}
    1 & \mbox{if} \ |x-m_i| \leq r_i \\
    0 & \mbox{otherwise}
    \end{array}
    \right.
\end{equation}
Using these,
we can construct a piecewise constant function $f:S^1\to \mathbb N$ by adding up the individual signal functions and restricting to the unit circle:
\begin{equation}\label{e-f-def-sum}
 f(x):=\sum_{i=1}^n h_i\big|_{S^1}(x)\;.
\end{equation}
This corresponds to the signal strength that an agent moving around the circle senses on its path.

Now, we want to consider the inverse problem: Given some piecewise constant $f:S^1\to \mathbb N$, reconstruct the data $(m_i,r_i)$, $i=1,\dots, n$, where the $m_i$'s and $r_i$'s as well as their number $n$ are a priori unknown. Furthermore, we want to do this based on incomplete information: we assume that we don't know the values of $f$ at a finite number of points. The reason for this formulation is to allow for uncertainty in the agent's measurements at the exact points where signal functions $h_i$ switch from 0 to 1. As we will see later, discontinuity points of each $h_i\big|_{S^1}$ will also be discontinuity points of (\ref{e-f-def-sum}), and specific choices of the values of the sum in~(\ref{e-f-def-sum}) at its discontinuity points will make it more desirable for our purposes.  Accordingly, we introduce the following terminology.

\begin{Definition}\label{d-decomposition}
 A sum of signal functions is called a \textit{decomposition} of some $f: S^1 \rightarrow \mathbb{N}$  if the sum of signal functions and $f$ disagree only on finitely many points in $S^1$.
\end{Definition}

This notion of decomposition gives rise to an equivalence relation on the set of piecewise constant functions from $S^1$ to $\mathbb{N}$; given a function $f$, we denote by $\bar f$ the equivalence class of all such functions which disagree with $f$ on a finite number of points in $S^1$. We will consider $f$ with finitely many discontinuities, so that we can find decompositions of the form (\ref{e-f-def-sum}) with finitely many signal functions.
Thus, given $\bar f$, we are tasked with finding all possible collections of data $(m_i,r_i)$, $i=1,\dots, n$ such that the corresponding sum~\eqref{e-f-def-sum} is a decomposition of $f$, i.e., is in $\bar f$. Furthermore, among such decompositions, we will be looking for those that satisfy an additional requirement of robustness and maximize some measure of likelihood, which we will formulate in the following sections.

We can also represent a piecewise constant ${f}: S^1 \rightarrow \mathbb{N}$ by the sequence of values it takes between discontinuity points, listed in the counterclockwise order (for concreteness) as the string $(a_1,...,a_{\ell})$.
To make this representation independent of the value we start with, we identify strings obtained from one another by a shift. This gives an equivalence class of strings,
\begin{equation}\label{stringclass}
    [a_1,\dots,a_{\ell}]=[a_{r\text{\,mod\,}\ell},\dots,a_{(r+\ell-1)\text{\,mod\,}\ell}]\quad \forall r\in \mathbb{Z}\;.
\end{equation}
Now, given a piecewise constant function on $S^1$ with $\ell$ discontinuity points, we define the map that outputs its equivalence class of strings as above,
\begin{equation}\label{seqmapdefn}
\text{seq}:= (f:S^1\rightarrow \mathbb{N})\mapsto [a_1,\dots,a_\ell]\;.
\end{equation}
We will call the image of $f$ under $\text{seq}(\cdot)$ the \textit{sequence representation} of $f$. Note that even if two functions belong to $\bar f$ (i.e., they disagree only on finitely many points) they need not have the same sequence representation. (As an example, consider the function $f\equiv 1$ with sequence representation $[1]$ and a function with two removable discontinuities and sequence representation $[1,1]$.) In other words, we have two different equivalence relations, one being the equivalence of strings and the other being the equivalence of functions.

It will be useful to have more explicit notation for the operator that builds the function $f$ from the individual signal functions as in~\eqref{e-f-def-sum}:
\begin{equation}\label{coupledsignal}
H_{n} :\mathbb{R}^{3n}\rightarrow
\{f:S^1\rightarrow \mathbb{N}\}, \quad
H_{n}(m_1,r_1,\dots,m_{n},r_{n}):=\sum_{i=1}^n  h_i\;.
\end{equation}
In what follows, when we write ``decomposition $H_n(v)$" we mean the sum of signal functions defined by~\eqref{e-h-binary-def} and~\eqref{coupledsignal} for some (not always explicitly specified) data vector $v=(m_1,r_1,\dots,m_{n},r_{n})$.
$H_{\{\cdot\}}$ is additive in the sense that
$
H_{m}(u) + H_{n}(v)= H_{(m+n)}(u,v).
$
In addition to using the standard symbol $\bdry$ for the boundary of a set, we denote by $\bdry h_i$ the intersection of the boundary of the {support} of a signal function $h_i$ with $S^1$, and we also extend this notation to a sum $H_n(v)$:
\begin{equation*}\label{bdry defintion}
\bdry h_i := \bdry(h_i^{-1}(1))\cap S^1,
 \quad  \bdry H_{n}(v):= \bigcup_{i=1}^n\bdry h_i\;.
\end{equation*}
We will call elements of the set $\bdry H_n(v)$ the \textit{boundary points} of $H_n(v)$. We use $\card$ to denote cardinality of a set (i.e., the number of distinct elements in the set).
Hence, for example, $\card \bdry h_i=2$ implies two distinct points of intersection, thus the support of such an $h_i$ intersected with $S^1$ is an arc of the circle (in other words, $h_i$ restricted to $S^1$ is identically 1 on some arc of the circle and 0 elsewhere). As all $h_i$ in (\ref{coupledsignal}) are compactly supported, a connection between the boundary points of $H_n(v)$ and the discontinuity points of $h_i$ is as follows.

\begin{Lemma}\label{l-bdry-points} ${x\in S^1}$ is a discontinuity point of $H_n(v)$ if and only if $x\in \bdry h_i$ for some $h_i$ with ${h_i^{-1}(1)\not\supseteq S^1}$.\\$\phantom \omitted$
\end{Lemma}

The proof of this lemma is elementary and is omitted. From now on, we will be using the symbol '$\triangle$' to indicate  results with omitted proof due to space limitations or because they are elementary.

\section{Robust decompositions}\label{s-robust}

The following notion of robustness will play a central role in the paper.

\begin{Definition}\label{d-robust}
A decomposition $H_n(v)$ is said to be \textit{robust} if sequence representations are invariant under {sufficiently small} perturbations, that is, if \[
\exists \bar\eps>0\; \text{\,s.t.}\;\forall \eps \in (0,\bar\eps) \;\forall u\in \mathbb{R}^{3n} \;\text{with}\;\lvert u\rvert=1,\;\text{seq} (H_{n}(v+\eps u)) = \text{seq}( H_{n}(v))\;.\]
\end{Definition}
This can be read as follows: $H_n(v)$ is robust when there exists an $\bar\eps>0$ such that for any perturbation of size less than $\bar\eps$ to the vector $v\in\R^{3n}$, the sequence representation of the perturbed decomposition is the same as the nominal decomposition.  In view of the definition~\eqref{seqmapdefn} of the map $\text{seq}(\cdot)$ and Lemma~\ref{l-bdry-points}, whether a decomposition is robust must depend on how the set of boundary points $\bdry H_n(v)$ changes under perturbations. To understand which decompositions are robust, we thus need to study how the sets $\bdry h_i$ behave under perturbations of $m_i$ and $r_i$.
Since $\bdry h_i$ is an intersection of two circles in the plane, it can have $0,1,2$ or infinitely many elements. Using elementary continuity arguments, we can easily see which values of $\card\bdry h_i$ lead to robustness.

\begin{Lemma}\label{l-element-counts-for-robustness}
The following statements hold:
\begin{itemize}
    \item[(i)]$\label{robustsignal}\text{If}\; \card \bdry h_i = 0 \ \text{or}\ 2,\:\text{then} $\[\exists \bar\eps>0\; \text{\,s.t.}\;\forall \eps \in (0,\bar\eps) \;\forall u\in \mathbb{R}^{3} \;\text{with}\;\lvert u\rvert=1 \;\text{we have}\; \card\bdry h_i(v+\eps u)=\card\bdry h_i(v)\]
    \item[(ii)]$\label{nonrobustsignal}\text{If}\; \card \bdry h_i=1\text{ or}\; \infty,\: \text{then for each }k\in\{0,1,2\},$ \[\forall \bar\eps>0 \;\exists \eps \in (0,\bar\eps)\;\exists u\; \text{with}\;\lvert u\rvert=1\; \text{\,s.t.}\;\card\bdry h_i(v+\eps u)=k\]\vspace{1ex}
\end{itemize}
\end{Lemma}

\begin{figure}[htbp]
\begin{tikzpicture}[x=0.60pt,y=0.60pt,yscale=-0.9,xscale=0.9]

\draw  [draw opacity=0][fill={rgb, 255:red, 34; green, 208; blue, 77 }  ,fill opacity=0.4 ][dash pattern={on 4.5pt off 4.5pt}] (34.33,217.88) .. controls (54.52,155.72) and (112.9,110.79) .. (181.79,110.79) .. controls (267.39,110.79) and (336.79,180.19) .. (336.79,265.79) .. controls (336.79,282.51) and (334.14,298.61) .. (329.24,313.7) -- (181.79,265.79) -- cycle ; \draw [color={rgb, 255:red, 98; green, 164; blue, 27 }  ,draw opacity=1 ][fill={rgb, 255:red, 34; green, 208; blue, 77 }  ,fill opacity=0.4 ][dash pattern={on 4.5pt off 4.5pt}] [dash pattern={on 4.5pt off 4.5pt}]  (34.33,217.88) .. controls (54.52,155.72) and (112.9,110.79) .. (181.79,110.79) .. controls (267.39,110.79) and (336.79,180.19) .. (336.79,265.79) .. controls (336.79,282.51) and (334.14,298.61) .. (329.24,313.7) ;
\draw  [draw opacity=0][fill={rgb, 255:red, 245; green, 166; blue, 35 }  ,fill opacity=0.4 ][dash pattern={on 5.63pt off 4.5pt}][line width=1.5]  (100.71,193.8) .. controls (123.05,163.52) and (159.11,143.86) .. (199.8,143.86) .. controls (267.59,143.86) and (322.55,198.45) .. (322.55,265.79) .. controls (322.55,292.72) and (313.76,317.6) .. (298.88,337.78) -- (199.8,265.79) -- cycle ; \draw [color={rgb, 255:red, 115; green, 104; blue, 67 }  ,draw opacity=1 ][fill={rgb, 255:red, 245; green, 166; blue, 35 }  ,fill opacity=0.4 ][dash pattern={on 5.63pt off 4.5pt}][line width=1.5]  [dash pattern={on 5.63pt off 4.5pt}]  (100.71,193.8) .. controls (123.05,163.52) and (159.11,143.86) .. (199.8,143.86) .. controls (267.59,143.86) and (322.55,198.45) .. (322.55,265.79) .. controls (322.55,292.72) and (313.76,317.6) .. (298.88,337.78) ;
\draw  [draw opacity=0][fill={rgb, 255:red, 78; green, 142; blue, 235 }  ,fill opacity=0.4 ][dash pattern={on 4.5pt off 4.5pt}] (158.59,193.8) .. controls (173.29,185.4) and (190.33,180.61) .. (208.5,180.61) .. controls (263.9,180.61) and (308.82,225.22) .. (308.82,280.25) .. controls (308.82,317.23) and (288.53,349.51) .. (258.41,366.7) -- (208.5,280.25) -- cycle ; \draw [color={rgb, 255:red, 56; green, 110; blue, 171 }  ,draw opacity=1 ][fill={rgb, 255:red, 78; green, 142; blue, 235 }  ,fill opacity=0.4 ][dash pattern={on 4.5pt off 4.5pt}] [dash pattern={on 4.5pt off 4.5pt}]  (158.59,193.8) .. controls (173.29,185.4) and (190.33,180.61) .. (208.5,180.61) .. controls (263.9,180.61) and (308.82,225.22) .. (308.82,280.25) .. controls (308.82,317.23) and (288.53,349.51) .. (258.41,366.7) ;
\draw  [color={rgb, 255:red, 150; green, 30; blue, 30 }  ,draw opacity=0.83 ][fill={rgb, 255:red, 229; green, 62; blue, 62 }  ,fill opacity=0.44 ][dash pattern={on 4.5pt off 4.5pt}] (181.79,265.79) .. controls (181.79,255.91) and (189.85,247.9) .. (199.8,247.9) .. controls (209.74,247.9) and (217.81,255.91) .. (217.81,265.79) .. controls (217.81,275.67) and (209.74,283.68) .. (199.8,283.68) .. controls (189.85,283.68) and (181.79,275.67) .. (181.79,265.79) -- cycle ;
\draw  [draw opacity=0] (392.83,205.36) .. controls (371.73,228.3) and (341.38,242.68) .. (307.64,242.68) .. controls (243.92,242.68) and (192.27,191.38) .. (192.27,128.09) .. controls (192.27,99.92) and (202.51,74.13) .. (219.48,54.17) -- (307.64,128.09) -- cycle ; \draw    (392.83,205.36) .. controls (371.73,228.3) and (341.38,242.68) .. (307.64,242.68) .. controls (243.92,242.68) and (192.27,191.38) .. (192.27,128.09) .. controls (192.27,99.92) and (202.51,74.13) .. (219.48,54.17) ;
\draw [fill={rgb, 255:red, 21; green, 19; blue, 19 }  ,fill opacity=1 ]   (199.8,265.79) -- (205.83,275.82) -- (207.47,278.54) ;
\draw [shift={(208.5,280.25)}, rotate = 238.96] [color={rgb, 255:red, 0; green, 0; blue, 0 }  ][line width=0.75]    (4.37,-1.32) .. controls (2.78,-0.56) and (1.32,-0.12) .. (0,0) .. controls (1.32,0.12) and (2.78,0.56) .. (4.37,1.32)   ;
\draw [shift={(199.8,265.79)}, rotate = 58.96] [color={rgb, 255:red, 0; green, 0; blue, 0 }  ][fill={rgb, 255:red, 0; green, 0; blue, 0 }  ][line width=0.75]      (0, 0) circle [x radius= 1.34, y radius= 1.34]   ;
\draw    (288.13,226.78) -- (331.87,205.22) ;
\draw [shift={(333.67,204.33)}, rotate = 153.76] [color={rgb, 255:red, 0; green, 0; blue, 0 }  ][line width=0.75]    (7.65,-3.43) .. controls (4.86,-1.61) and (2.31,-0.47) .. (0,0) .. controls (2.31,0.47) and (4.86,1.61) .. (7.65,3.43)   ;
\draw [shift={(286.33,227.67)}, rotate = 333.76] [color={rgb, 255:red, 0; green, 0; blue, 0 }  ][line width=0.75]    (7.65,-3.43) .. controls (4.86,-1.61) and (2.31,-0.47) .. (0,0) .. controls (2.31,0.47) and (4.86,1.61) .. (7.65,3.43)   ;

\draw (190.83,270.4) node [anchor=north west][inner sep=0.75pt]  [font=\footnotesize]  {$\varepsilon _{1}$};
\draw (317,211.4) node [anchor=north west][inner sep=0.75pt]  [font=\small]  {$\varepsilon _{2}$};
\draw (214,70.07) node [anchor=north west][inner sep=0.75pt]  [font=\large]  {$S^{1}$};
\draw (55.83,208.66) node [anchor=north west][inner sep=0.75pt]  [font=\normalsize]  {$\#\partial h_{i} =2$};

\draw [color={rgb, 255:red, 65; green, 117; blue, 5 }  ,draw opacity=1 ][fill={rgb, 255:red, 227; green, 227; blue, 48 }  ,fill opacity=1 ]  (334.57, 239.54) circle [x radius= 3, y radius= 3]   ;
\draw [color={rgb, 255:red, 65; green, 117; blue, 5 }  ,draw opacity=1 ][fill={rgb, 255:red, 227; green, 227; blue, 48 }  ,fill opacity=1 ]  (193.51, 111.23) circle [x radius= 3, y radius= 3]   ;
\draw [color={rgb, 255:red, 65; green, 117; blue, 5 }  ,draw opacity=1 ][fill={rgb, 255:red, 227; green, 227; blue, 48 }  ,fill opacity=1 ]  (320.22, 242.01) circle [x radius= 3, y radius= 3]   ;
\draw [color={rgb, 255:red, 65; green, 117; blue, 5 }  ,draw opacity=1 ][fill={rgb, 255:red, 227; green, 227; blue, 48 }  ,fill opacity=1 ]  (193.38, 144.02) circle [x radius= 3, y radius= 3]   ;
\draw [color={rgb, 255:red, 65; green, 117; blue, 5 }  ,draw opacity=1 ][fill={rgb, 255:red, 227; green, 227; blue, 48 }  ,fill opacity=1 ]  (205.1, 180.66) circle [x radius= 3, y radius= 3]   ;
\draw [color={rgb, 255:red, 65; green, 117; blue, 5 }  ,draw opacity=1 ][fill={rgb, 255:red, 227; green, 227; blue, 48 }  ,fill opacity=1 ]  (301.37, 242.52) circle [x radius= 3, y radius= 3]   ;
\end{tikzpicture}
\hspace{10ex}
\begin{tikzpicture}[x=0.58pt,y=0.58pt,yscale=-0.9,xscale=0.9]

\draw  [draw opacity=0][fill={rgb, 255:red, 34; green, 208; blue, 77 }  ,fill opacity=0.4 ][dash pattern={on 4.5pt off 4.5pt}] (64.34,197.72) .. controls (88.94,144.07) and (146.02,109.85) .. (207.82,117.05) .. controls (284.62,125.99) and (339.63,195.49) .. (330.69,272.29) .. controls (328.95,287.3) and (324.89,301.47) .. (318.92,314.49) -- (191.63,256.11) -- cycle ; \draw [color={rgb, 255:red, 98; green, 164; blue, 27 }  ,draw opacity=1 ][fill={rgb, 255:red, 34; green, 208; blue, 77 }  ,fill opacity=0.4 ][dash pattern={on 4.5pt off 4.5pt}] [dash pattern={on 4.5pt off 4.5pt}]  (64.34,197.72) .. controls (88.94,144.07) and (146.02,109.85) .. (207.82,117.05) .. controls (284.62,125.99) and (339.63,195.49) .. (330.69,272.29) .. controls (328.95,287.3) and (324.89,301.47) .. (318.92,314.49) ;
\draw  [draw opacity=0][fill={rgb, 255:red, 245; green, 166; blue, 35 }  ,fill opacity=0.4 ][dash pattern={on 5.63pt off 4.5pt}][line width=1.5]  (113.67,174.45) .. controls (137.72,155.05) and (170.39,147.02) .. (202.51,155.66) .. controls (255.84,170) and (287.45,224.86) .. (273.11,278.19) .. controls (267.41,299.41) and (255.29,317.19) .. (239.41,330) -- (176.54,252.23) -- cycle ; \draw [color={rgb, 255:red, 115; green, 104; blue, 67 }  ,draw opacity=1 ][fill={rgb, 255:red, 245; green, 166; blue, 35 }  ,fill opacity=0.4 ][dash pattern={on 5.63pt off 4.5pt}][line width=1.5]  [dash pattern={on 5.63pt off 4.5pt}]  (113.67,174.45) .. controls (137.72,155.05) and (170.39,147.02) .. (202.51,155.66) .. controls (255.84,170) and (287.45,224.86) .. (273.11,278.19) .. controls (267.41,299.41) and (255.29,317.19) .. (239.41,330) ;
\draw  [draw opacity=0] (412.58,206.86) .. controls (387.73,225.9) and (355.25,234.79) .. (322.02,228.76) .. controls (259.53,217.4) and (218.08,157.54) .. (229.44,95.05) .. controls (234.47,67.37) and (249.02,43.81) .. (269.09,27.15) -- (342.58,115.61) -- cycle ; \draw    (412.58,206.86) .. controls (387.73,225.9) and (355.25,234.79) .. (322.02,228.76) .. controls (259.53,217.4) and (218.08,157.54) .. (229.44,95.05) .. controls (234.47,67.37) and (249.02,43.81) .. (269.09,27.15) ;
\draw [fill={rgb, 255:red, 248; green, 231; blue, 28 }  ,fill opacity=1 ]   (175.8,248) -- (181.83,258.03) -- (183.47,260.75) ;
\draw [shift={(184.5,262.46)}, rotate = 238.96] [color={rgb, 255:red, 0; green, 0; blue, 0 }  ][line width=0.75]    (4.37,-1.32) .. controls (2.78,-0.56) and (1.32,-0.12) .. (0,0) .. controls (1.32,0.12) and (2.78,0.56) .. (4.37,1.32)   ;
\draw [shift={(175.8,248)}, rotate = 58.96] [color={rgb, 255:red, 0; green, 0; blue, 0 }  ][fill={rgb, 255:red, 0; green, 0; blue, 0 }  ][line width=0.75]      (0, 0) circle [x radius= 1.34, y radius= 1.34]   ;
\draw  [draw opacity=0][fill={rgb, 255:red, 74; green, 144; blue, 226 }  ,fill opacity=0.4 ][dash pattern={on 4.5pt off 4.5pt}] (145.95,170.91) .. controls (183.41,164.87) and (222.09,183.13) .. (240.35,218.72) .. controls (263.03,262.95) and (245.56,317.19) .. (201.33,339.87) .. controls (192.69,344.31) and (183.67,347.2) .. (174.58,348.67) -- (160.26,259.79) -- cycle ; \draw [color={rgb, 255:red, 74; green, 144; blue, 226 }  ,draw opacity=1 ][fill={rgb, 255:red, 74; green, 144; blue, 226 }  ,fill opacity=0.4 ][dash pattern={on 4.5pt off 4.5pt}] [dash pattern={on 4.5pt off 4.5pt}]  (145.95,170.91) .. controls (183.41,164.87) and (222.09,183.13) .. (240.35,218.72) .. controls (263.03,262.95) and (245.56,317.19) .. (201.33,339.87) .. controls (192.69,344.31) and (183.67,347.2) .. (174.58,348.67) ;
\draw  [color={rgb, 255:red, 150; green, 30; blue, 30 }  ,draw opacity=0.83 ][fill={rgb, 255:red, 229; green, 62; blue, 62 }  ,fill opacity=0.44 ][dash pattern={on 4.5pt off 4.5pt}] (157.79,248) .. controls (157.79,238.12) and (165.85,230.11) .. (175.8,230.11) .. controls (185.74,230.11) and (193.81,238.12) .. (193.81,248) .. controls (193.81,257.88) and (185.74,265.89) .. (175.8,265.89) .. controls (165.85,265.89) and (157.79,257.88) .. (157.79,248) -- cycle ;
\draw  [fill={rgb, 255:red, 248; green, 231; blue, 28 }  ,fill opacity=1 ] (249,186.5) .. controls (249,184.9) and (250.3,183.6) .. (251.9,183.6) .. controls (253.5,183.6) and (254.8,184.9) .. (254.8,186.5) .. controls (254.8,188.1) and (253.5,189.4) .. (251.9,189.4) .. controls (250.3,189.4) and (249,188.1) .. (249,186.5) -- cycle ;
\draw    (245.35,125.07) -- (210.65,186.26) ;
\draw [shift={(209.67,188)}, rotate = 299.55] [color={rgb, 255:red, 0; green, 0; blue, 0 }  ][line width=0.75]    (10.93,-3.29) .. controls (6.95,-1.4) and (3.31,-0.3) .. (0,0) .. controls (3.31,0.3) and (6.95,1.4) .. (10.93,3.29)   ;
\draw [shift={(246.33,123.33)}, rotate = 119.55] [color={rgb, 255:red, 0; green, 0; blue, 0 }  ][line width=0.75]    (10.93,-3.29) .. controls (6.95,-1.4) and (3.31,-0.3) .. (0,0) .. controls (3.31,0.3) and (6.95,1.4) .. (10.93,3.29)   ;

\draw (166.83,252.61) node [anchor=north west][inner sep=0.75pt]  [font=\footnotesize]  {$\varepsilon _{1}$};
\draw (254,59.28) node [anchor=north west][inner sep=0.75pt]  [font=\large]  {$S^{1}$};
\draw (177.62,351.18) node [anchor=north west][inner sep=0.75pt]  [font=\scriptsize,color={rgb, 255:red, 39; green, 82; blue, 132 }  ,opacity=1 ,rotate=-340]  {$\#\partial h_{i} =0$};
\draw (242.45,332.51) node [anchor=north west][inner sep=0.75pt]  [font=\scriptsize,color={rgb, 255:red, 71; green, 66; blue, 11 }  ,opacity=1 ,rotate=-340]  {$\#\partial h_{i} =1$};
\draw (321.96,317) node [anchor=north west][inner sep=0.75pt]  [font=\scriptsize,color={rgb, 255:red, 58; green, 98; blue, 10 }  ,opacity=1 ,rotate=-340]  {$\#\partial h_{i} =2$};
\draw (237.67,136.73) node [anchor=north west][inner sep=0.75pt]  [font=\small]  {$\varepsilon _{2}$};

\draw [fill={rgb, 255:red, 248; green, 231; blue, 28 }  ,fill opacity=1 ]  (329.18, 229.83) circle [x radius= 3, y radius= 3]   ;
\draw [fill={rgb, 255:red, 248; green, 231; blue, 28 }  ,fill opacity=1 ]  (227.68, 120.81) circle [x radius= 3, y radius= 3]   ;
\end{tikzpicture}
\caption{\small Examples with $\card\bdry h_i=2$ (left) and $\card\bdry h_i=1$ (right) demonstrating small perturbations of $h_i$. The dashed arcs represent $\bdry(h_i^{-1}(1))$ in the nominal case.}
\label{fig:perturbations}
\end{figure}
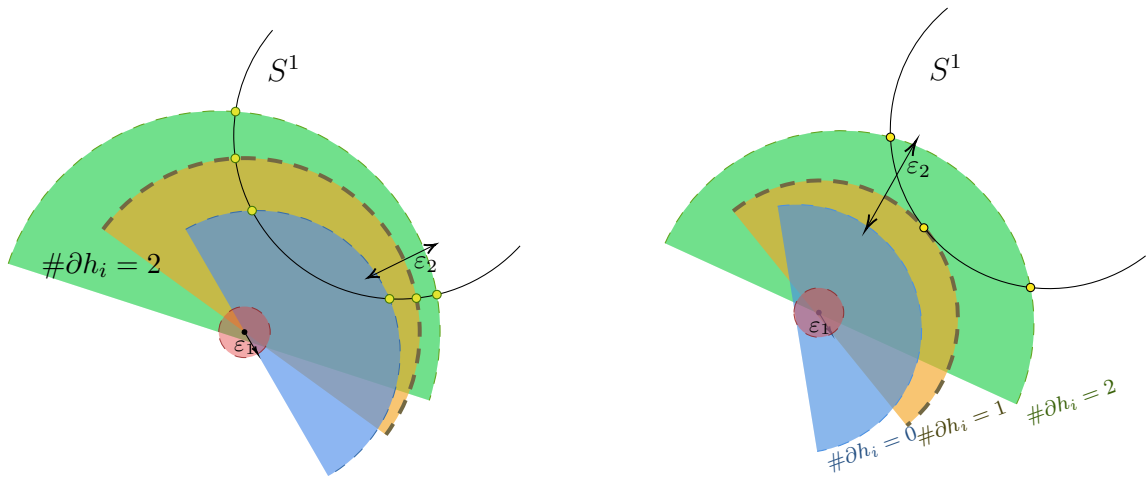

The proof of Lemma~\ref{l-element-counts-for-robustness} is straightforward and we omit it; instead, we give Figure~\ref{fig:perturbations} as illustration.

A change in $\card \bdry h_i$ corresponds to a change in the sequence representation, thus causing non-robustness. By statement~(ii) of Lemma~\ref{l-element-counts-for-robustness}, signal functions whose support's boundary intersects $S^1$ at one or infinitely many points are non-robust decompositions by themselves. Conversely, by statement~(i) of the lemma, for $\card\bdry h_i\in \{0,2\}$ there always exists an $\bar\eps>0$ such that if $h_i$ is perturbed less that $\bar\eps$, $\card \bdry h_i$ is preserved. Therefore, $h_i$ with $\card \bdry h_i \in \{0,2\}$  are robust decompositions by themselves. Considering this fact, we will refer to $h_i$ with $\card \bdry h_i \in \{0,2\} $ as \textit{robust signals}. However, even if we use only robust signals in (\ref{coupledsignal}), we still cannot expect robustness of $H_n(v)$, as shown in the following example.

\begin{Example}\label{x-non-robust-two-cases} Consider the case when $n = 2$ and $\card\bdry H_2$ is an odd number. The following specific cases (and subcases) are possible:
\begin{enumerate}
  \item $\card\bdry H_2(v)=1$, and either
  $\card\partial h_1=1, \card\partial h_2=0\; \text{(or vice versa)}$ or
   $\card\partial h_1=\card\partial h_2=1$;
  \item $\card\bdry H_2(v)=3$, and either
 $\card\partial h_1=1, \card\partial h_2=2$\;\text{(or vice versa)} or
    $\card\partial h_1=\card\partial h_2=2$.

\end{enumerate}
In all of the above cases except the very last subcase, we have at least one signal with $\card\bdry h_i=1$. Such an $H_2(v)$ cannot be robust because for perturbations only affecting the radius $r_i$ of one of the signals with $\card \bdry h_i=1$, the sequence representation will change. In the last subcase, it is clear that $\bdry h_1 \cap \bdry h_2\neq \emptyset$ is inevitable. To see why the presence of such shared boundary points causes non-robustness, consider Figure~\ref{f-shared-bdry-perturb}.
Similar considerations show that, more generally, $H_n(v)$ cannot be robust for odd values of $\card\bdry H_n(v)$; we will be able to provide a proof of this fact from the subsequent discussion, as a corollary of Proposition 3.2.

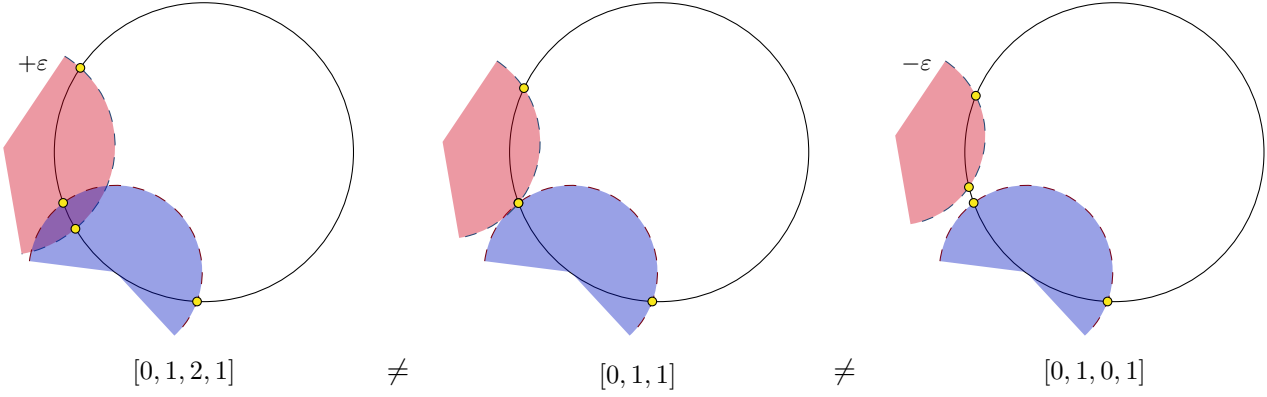
\begin{figure}
\begin{tikzpicture}[x=0.82pt,y=0.82pt,yscale=-1,xscale=1]

\draw   (297.05,129) .. controls (297.05,90.89) and (327.94,60) .. (366.05,60) .. controls (404.16,60) and (435.05,90.89) .. (435.05,129) .. controls (435.05,167.11) and (404.16,198) .. (366.05,198) .. controls (327.94,198) and (297.05,167.11) .. (297.05,129) -- cycle ;
\draw  [draw opacity=0][fill={rgb, 255:red, 2; green, 20; blue, 193 }  ,fill opacity=0.4 ][dash pattern={on 4.5pt off 4.5pt}] (285.57,179.3) .. controls (286.14,174.79) and (287.49,170.31) .. (289.69,166.03) .. controls (299.79,146.39) and (323.91,138.65) .. (343.55,148.75) .. controls (363.2,158.85) and (370.94,182.97) .. (360.84,202.61) .. controls (358.64,206.89) and (355.78,210.6) .. (352.44,213.68) -- (325.26,184.32) -- cycle ; \draw  [color={rgb, 255:red, 120; green, 2; blue, 17 }  ,draw opacity=1 ][dash pattern={on 4.5pt off 4.5pt}] (285.57,179.3) .. controls (286.14,174.79) and (287.49,170.31) .. (289.69,166.03) .. controls (299.79,146.39) and (323.91,138.65) .. (343.55,148.75) .. controls (363.2,158.85) and (370.94,182.97) .. (360.84,202.61) .. controls (358.64,206.89) and (355.78,210.6) .. (352.44,213.68) ;
\draw  [draw opacity=0][fill={rgb, 255:red, 208; green, 2; blue, 27 }  ,fill opacity=0.4 ][dash pattern={on 4.5pt off 4.5pt}] (291.21,86.87) .. controls (299.69,92.56) and (306.3,101.15) .. (309.36,111.73) .. controls (316.27,135.61) and (302.52,160.56) .. (278.65,167.47) .. controls (277.03,167.94) and (275.41,168.31) .. (273.8,168.59) -- (266.14,124.24) -- cycle ; \draw  [color={rgb, 255:red, 38; green, 74; blue, 118 }  ,draw opacity=1 ][dash pattern={on 4.5pt off 4.5pt}] (291.21,86.87) .. controls (299.69,92.56) and (306.3,101.15) .. (309.36,111.73) .. controls (316.27,135.61) and (302.52,160.56) .. (278.65,167.47) .. controls (277.03,167.94) and (275.41,168.31) .. (273.8,168.59) ;
\draw   (87.05,129) .. controls (87.05,90.89) and (117.94,60) .. (156.05,60) .. controls (194.16,60) and (225.05,90.89) .. (225.05,129) .. controls (225.05,167.11) and (194.16,198) .. (156.05,198) .. controls (117.94,198) and (87.05,167.11) .. (87.05,129) -- cycle ;
\draw  [draw opacity=0][fill={rgb, 255:red, 208; green, 2; blue, 27 }  ,fill opacity=0.4 ][dash pattern={on 4.5pt off 4.5pt}] (91.84,84.69) .. controls (101.9,90.88) and (109.69,100.53) .. (113.17,112.55) .. controls (120.81,138.94) and (104.76,166.77) .. (77.33,174.71) .. controls (75.54,175.22) and (73.76,175.64) .. (71.97,175.97) -- (63.5,126.93) -- cycle ; \draw  [color={rgb, 255:red, 38; green, 74; blue, 118 }  ,draw opacity=1 ][dash pattern={on 4.5pt off 4.5pt}] (91.84,84.69) .. controls (101.9,90.88) and (109.69,100.53) .. (113.17,112.55) .. controls (120.81,138.94) and (104.76,166.77) .. (77.33,174.71) .. controls (75.54,175.22) and (73.76,175.64) .. (71.97,175.97) ;
\draw   (507.05,129) .. controls (507.05,90.89) and (537.94,60) .. (576.05,60) .. controls (614.16,60) and (645.05,90.89) .. (645.05,129) .. controls (645.05,167.11) and (614.16,198) .. (576.05,198) .. controls (537.94,198) and (507.05,167.11) .. (507.05,129) -- cycle ;
\draw  [draw opacity=0][fill={rgb, 255:red, 2; green, 20; blue, 193 }  ,fill opacity=0.4 ][dash pattern={on 4.5pt off 4.5pt}] (495.57,179.3) .. controls (496.14,174.79) and (497.49,170.31) .. (499.69,166.03) .. controls (509.79,146.39) and (533.91,138.65) .. (553.55,148.75) .. controls (573.2,158.85) and (580.94,182.97) .. (570.84,202.61) .. controls (568.64,206.89) and (565.78,210.6) .. (562.44,213.68) -- (535.26,184.32) -- cycle ; \draw  [color={rgb, 255:red, 120; green, 2; blue, 27 }  ,draw opacity=1 ][dash pattern={on 4.5pt off 4.5pt}] (495.57,179.3) .. controls (496.14,174.79) and (497.49,170.31) .. (499.69,166.03) .. controls (509.79,146.39) and (533.91,138.65) .. (553.55,148.75) .. controls (573.2,158.85) and (580.94,182.97) .. (570.84,202.61) .. controls (568.64,206.89) and (565.78,210.6) .. (562.44,213.68) ;
\draw  [draw opacity=0][fill={rgb, 255:red, 208; green, 2; blue, 27 }  ,fill opacity=0.4 ][dash pattern={on 4.5pt off 4.5pt}] (498.06,86.88) .. controls (505.88,92.12) and (511.97,100.05) .. (514.8,109.8) .. controls (521.17,131.82) and (508.49,154.83) .. (486.47,161.2) .. controls (484.98,161.63) and (483.49,161.98) .. (482,162.24) -- (474.93,121.34) -- cycle ; \draw  [color={rgb, 255:red, 38; green, 74; blue, 118 }  ,draw opacity=1 ][dash pattern={on 4.5pt off 4.5pt}] (498.06,86.88) .. controls (505.88,92.12) and (511.97,100.05) .. (514.8,109.8) .. controls (521.17,131.82) and (508.49,154.83) .. (486.47,161.2) .. controls (484.98,161.63) and (483.49,161.98) .. (482,162.24) ;
\draw  [draw opacity=0][fill={rgb, 255:red, 2; green, 20; blue, 193 }  ,fill opacity=0.4 ][dash pattern={on 4.5pt off 4.5pt}] (75.57,179.3) .. controls (76.14,174.79) and (77.49,170.31) .. (79.69,166.03) .. controls (89.79,146.39) and (113.91,138.65) .. (133.55,148.75) .. controls (153.2,158.85) and (160.94,182.97) .. (150.84,202.61) .. controls (148.64,206.89) and (145.78,210.6) .. (142.44,213.68) -- (115.26,184.32) -- cycle ; \draw  [color={rgb, 255:red, 120; green, 4; blue, 18 }  ,draw opacity=1 ][dash pattern={on 4.5pt off 4.5pt}] (75.57,179.3) .. controls (76.14,174.79) and (77.49,170.31) .. (79.69,166.03) .. controls (89.79,146.39) and (113.91,138.65) .. (133.55,148.75) .. controls (153.2,158.85) and (160.94,182.97) .. (150.84,202.61) .. controls (148.64,206.89) and (145.78,210.6) .. (142.44,213.68) ;

\draw (122,223.4) node [anchor=north west][inner sep=0.75pt]  [font=\small]  {$[ 0,1,2,1]$};
\draw (337,225.4) node [anchor=north west][inner sep=0.75pt]  [font=\small]  {$[ 0,1,1]$};
\draw (542,223.4) node [anchor=north west][inner sep=0.75pt]  [font=\small]  {$[ 0,1,0,1]$};
\draw (239.33,223.4) node [anchor=north west][inner sep=0.75pt]    {$\neq $};
\draw (445.33,223.4) node [anchor=north west][inner sep=0.75pt]    {$\neq $};
\draw (69,83.4) node [anchor=north west][inner sep=0.75pt]  [font=\footnotesize]  {$+\varepsilon $};
\draw (476.67,83.4) node [anchor=north west][inner sep=0.75pt]  [font=\footnotesize]  {$-\varepsilon $};

\draw [fill={rgb, 255:red, 248; green, 231; blue, 28 }  ,fill opacity=1 ]  (301.12, 152.41) circle [x radius= 2, y radius= 2]   ;
\draw [fill={rgb, 255:red, 248; green, 231; blue, 28 }  ,fill opacity=1 ]  (362.88, 197.93) circle [x radius= 2, y radius= 2]   ;
\draw [fill={rgb, 255:red, 248; green, 231; blue, 28 }  ,fill opacity=1 ]  (303.69, 99.43) circle [x radius= 2, y radius= 2]   ;
\draw [fill={rgb, 255:red, 248; green, 231; blue, 28 }  ,fill opacity=1 ]  (301.16, 152.51) circle [x radius= 2, y radius= 2]   ;
\draw [fill={rgb, 255:red, 248; green, 231; blue, 28 }  ,fill opacity=1 ]  (99.06, 90.08) circle [x radius= 2, y radius= 2]   ;
\draw [fill={rgb, 255:red, 248; green, 231; blue, 28 }  ,fill opacity=1 ]  (96.79, 164.37) circle [x radius= 2, y radius= 2]   ;
\draw [fill={rgb, 255:red, 248; green, 231; blue, 28 }  ,fill opacity=1 ]  (91.12, 152.41) circle [x radius= 2, y radius= 2]   ;
\draw [fill={rgb, 255:red, 248; green, 231; blue, 28 }  ,fill opacity=1 ]  (152.88, 197.93) circle [x radius= 2, y radius= 2]   ;
\draw [fill={rgb, 255:red, 248; green, 231; blue, 28 }  ,fill opacity=1 ]  (511.12, 152.41) circle [x radius= 2, y radius= 2]   ;
\draw [fill={rgb, 255:red, 248; green, 231; blue, 28 }  ,fill opacity=1 ]  (572.88, 197.93) circle [x radius= 2, y radius= 2]   ;
\draw [fill={rgb, 255:red, 248; green, 231; blue, 28 }  ,fill opacity=1 ]  (512.14, 102.94) circle [x radius= 2, y radius= 2]   ;
\draw [fill={rgb, 255:red, 248; green, 231; blue, 28 }  ,fill opacity=1 ]  (508.94, 145.12) circle [x radius= 2, y radius= 2]   ;
\end{tikzpicture}
\caption{\small A decomposition formed by two signals with supports intersecting at one point (center) and its perturbations (left and right).}\label{f-shared-bdry-perturb}
\end{figure}
\end{Example}
\begin{Example} Consider  $k h_i$, $k\in \mathbb{N}$ (i.e., the sum of $k$ identical signal functions). If $k=1$, then $h_i$ is robust if and only if $\bdry h_i$ has 0 or 2 elements (Lemma~\ref{l-element-counts-for-robustness}). If $k > 1$ then, unless $\card\bdry h_i=0$, perturbing the radius of just one signal function will always change $\card \bdry H_{k}(\cdot)$, thus changing the sequence representation; see Figure~\ref{fig:multiplicity}. Therefore, if $k >1$, $k h_i$ is robust if and only if $\card \bdry h_i=0$.

\begin{figure}[h!]
\begin{tikzpicture}[x=0.82pt,y=0.82pt,yscale=-1,xscale=1]

\draw  [draw opacity=0][fill={rgb, 255:red, 200; green, 0; blue, 0 }  ,fill opacity=0.4 ][dash pattern={on 4.5pt off 4.5pt}] (267.81,126.44) .. controls (296.34,119.79) and (327.36,130.3) .. (345.62,155.51) .. controls (363.9,180.77) and (364.17,213.59) .. (348.86,238.63) -- (284.87,199.49) -- cycle ; \draw  [color={rgb, 255:red, 120; green, 4; blue, 18 }  ,draw opacity=1 ][dash pattern={on 4.5pt off 4.5pt}] (267.81,126.44) .. controls (296.34,119.79) and (327.36,130.3) .. (345.62,155.51) .. controls (363.9,180.77) and (364.17,213.59) .. (348.86,238.63) ;
\draw  [draw opacity=0][fill={rgb, 255:red, 200; green, 2; blue, 0 }  ,fill opacity=0.4 ][dash pattern={on 4.5pt off 4.5pt}] (471.07,118.01) .. controls (505.3,110.03) and (542.53,122.64) .. (564.44,152.9) .. controls (586.38,183.21) and (586.69,222.6) .. (568.33,252.65) -- (491.53,205.68) -- cycle ; \draw  [color={rgb, 255:red, 208; green, 2; blue, 27 }  ,draw opacity=1 ][dash pattern={on 4.5pt off 4.5pt}] (471.07,118.01) .. controls (505.3,110.03) and (542.53,122.64) .. (564.44,152.9) .. controls (586.38,183.21) and (586.69,222.6) .. (568.33,252.65) ;
\draw  [draw opacity=0][fill={rgb, 255:red, 0; green, 200; blue, 0 }  ,fill opacity=0.4 ][dash pattern={on 4.5pt off 4.5pt}] (474.48,132.62) .. controls (503.01,125.97) and (534.03,136.48) .. (552.29,161.7) .. controls (570.57,186.96) and (570.83,219.78) .. (555.53,244.82) -- (491.53,205.68) -- cycle ; \draw  [color={rgb, 255:red, 65; green, 117; blue, 5 }  ,draw opacity=1 ][dash pattern={on 4.5pt off 4.5pt}] (474.48,132.62) .. controls (503.01,125.97) and (534.03,136.48) .. (552.29,161.7) .. controls (570.57,186.96) and (570.83,219.78) .. (555.53,244.82) ;
\draw   (289.05,149) .. controls (289.05,110.89) and (319.94,80) .. (358.05,80) .. controls (396.16,80) and (427.05,110.89) .. (427.05,149) .. controls (427.05,187.11) and (396.16,218) .. (358.05,218) .. controls (319.94,218) and (289.05,187.11) .. (289.05,149) -- cycle ;
\draw  [draw opacity=0][fill={rgb, 255:red, 0; green, 200; blue, 0 }  ,fill opacity=0.4 ][dash pattern={on 4.5pt off 4.5pt}] (267.81,126.44) .. controls (296.34,119.79) and (327.36,130.3) .. (345.62,155.51) .. controls (363.9,180.77) and (364.17,213.59) .. (348.86,238.63) -- (284.87,199.49) -- cycle ; \draw  [color={rgb, 255:red, 120; green, 4; blue, 18 }  ,draw opacity=1 ][dash pattern={on 4.5pt off 4.5pt}] (267.81,126.44) .. controls (296.34,119.79) and (327.36,130.3) .. (345.62,155.51) .. controls (363.9,180.77) and (364.17,213.59) .. (348.86,238.63) ;
\draw   (495.72,155.19) .. controls (495.72,117.08) and (526.61,86.19) .. (564.72,86.19) .. controls (602.82,86.19) and (633.72,117.08) .. (633.72,155.19) .. controls (633.72,193.29) and (602.82,224.19) .. (564.72,224.19) .. controls (526.61,224.19) and (495.72,193.29) .. (495.72,155.19) -- cycle ;
\draw  [draw opacity=0][fill={rgb, 255:red, 200; green, 0; blue, 0 }  ,fill opacity=0.4 ][dash pattern={on 4.5pt off 4.5pt}] (63.89,140.38) .. controls (86.71,135.06) and (111.53,143.47) .. (126.14,163.64) .. controls (140.76,183.85) and (140.97,210.1) .. (128.73,230.14) -- (77.53,198.82) -- cycle ; \draw  [color={rgb, 255:red, 255; green, 2; blue, 27 }  ,draw opacity=1 ][dash pattern={on 4.5pt off 4.5pt}] (63.89,140.38) .. controls (86.71,135.06) and (111.53,143.47) .. (126.14,163.64) .. controls (140.76,183.85) and (140.97,210.1) .. (128.73,230.14) ;
\draw   (81.72,149) .. controls (81.72,110.89) and (112.61,80) .. (150.72,80) .. controls (188.82,80) and (219.72,110.89) .. (219.72,149) .. controls (219.72,187.11) and (188.82,218) .. (150.72,218) .. controls (112.61,218) and (81.72,187.11) .. (81.72,149) -- cycle ;
\draw  [draw opacity=0][fill={rgb, 255:red, 0; green, 200; blue, 0 }  ,fill opacity=0.4 ][dash pattern={on 4.5pt off 4.5pt}] (60.48,125.77) .. controls (89.01,119.12) and (120.03,129.63) .. (138.29,154.85) .. controls (156.57,180.11) and (156.83,212.92) .. (141.53,237.97) -- (77.53,198.82) -- cycle ; \draw  [color={rgb, 255:red, 65; green, 117; blue, 5 }  ,draw opacity=1 ][dash pattern={on 4.5pt off 4.5pt}] (60.48,125.77) .. controls (89.01,119.12) and (120.03,129.63) .. (138.29,154.85) .. controls (156.57,180.11) and (156.83,212.92) .. (141.53,237.97) ;

\draw (42.67,122.73) node [anchor=north west][inner sep=0.75pt]  [font=\footnotesize]  {$h_{2}$};
\draw (35.33,143.4) node [anchor=north west][inner sep=0.75pt]  [font=\footnotesize]  {$h_{1}^{-\varepsilon }$};
\draw (240,122.07) node [anchor=north west][inner sep=0.75pt]  [font=\footnotesize]  {$h_{1,2}$};
\draw (458.67,133.4) node [anchor=north west][inner sep=0.75pt]  [font=\footnotesize]  {$h_{2}$};
\draw (447.33,113.4) node [anchor=north west][inner sep=0.75pt]  [font=\footnotesize]  {$h_{1}^{\ +\varepsilon }$};
\draw (112,253.4) node [anchor=north west][inner sep=0.75pt]  [font=\small]  {$[ 0,1,2,1]$};
\draw (542,253.4) node [anchor=north west][inner sep=0.75pt]  [font=\small]  {$[ 0,1,2,1]$};
\draw (342,253.73) node [anchor=north west][inner sep=0.75pt]  [font=\small]  {$[ 0,2]$};
\draw (237,253.4) node [anchor=north west][inner sep=0.75pt]    {$\neq $};
\draw (438.67,253.4) node [anchor=north west][inner sep=0.75pt]    {$\neq $};

\draw [fill={rgb, 255:red, 248; green, 231; blue, 28 }  ,fill opacity=1 ]  (293.35, 124.97) circle [x radius= 2, y radius= 2]   ;
\draw [fill={rgb, 255:red, 248; green, 231; blue, 28 }  ,fill opacity=1 ]  (357.55, 218) circle [x radius= 2, y radius= 2]   ;
\draw [fill={rgb, 255:red, 248; green, 231; blue, 28 }  ,fill opacity=1 ]  (293.35, 124.97) circle [x radius= 2, y radius= 2]   ;
\draw [fill={rgb, 255:red, 248; green, 231; blue, 28 }  ,fill opacity=1 ]  (357.55, 218) circle [x radius= 2, y radius= 2]   ;
\draw [fill={rgb, 255:red, 248; green, 231; blue, 28 }  ,fill opacity=1 ]  (507.2, 117.05) circle [x radius= 2, y radius= 2]   ;
\draw [fill={rgb, 255:red, 248; green, 231; blue, 28 }  ,fill opacity=1 ]  (579.95, 222.5) circle [x radius= 2, y radius= 2]   ;
\draw [fill={rgb, 255:red, 248; green, 231; blue, 28 }  ,fill opacity=1 ]  (500.02, 131.16) circle [x radius= 2, y radius= 2]   ;
\draw [fill={rgb, 255:red, 248; green, 231; blue, 28 }  ,fill opacity=1 ]  (564.22, 224.18) circle [x radius= 2, y radius= 2]   ;
\draw [fill={rgb, 255:red, 248; green, 231; blue, 28 }  ,fill opacity=1 ]  (82.43, 139.02) circle [x radius= 2, y radius= 2]   ;
\draw [fill={rgb, 255:red, 248; green, 231; blue, 28 }  ,fill opacity=1 ]  (134.97, 216.19) circle [x radius= 2, y radius= 2]   ;
\draw [fill={rgb, 255:red, 248; green, 231; blue, 28 }  ,fill opacity=1 ]  (86.26, 124.33) circle [x radius= 2, y radius= 2]   ;
\draw [fill={rgb, 255:red, 248; green, 231; blue, 28 }  ,fill opacity=1 ]  (150.04, 218) circle [x radius= 2, y radius= 2]   ;
\end{tikzpicture}
\caption{\small A decomposition using two identical signals (center) and its perturbations (left and right).}\label{fig:multiplicity}
\end{figure}
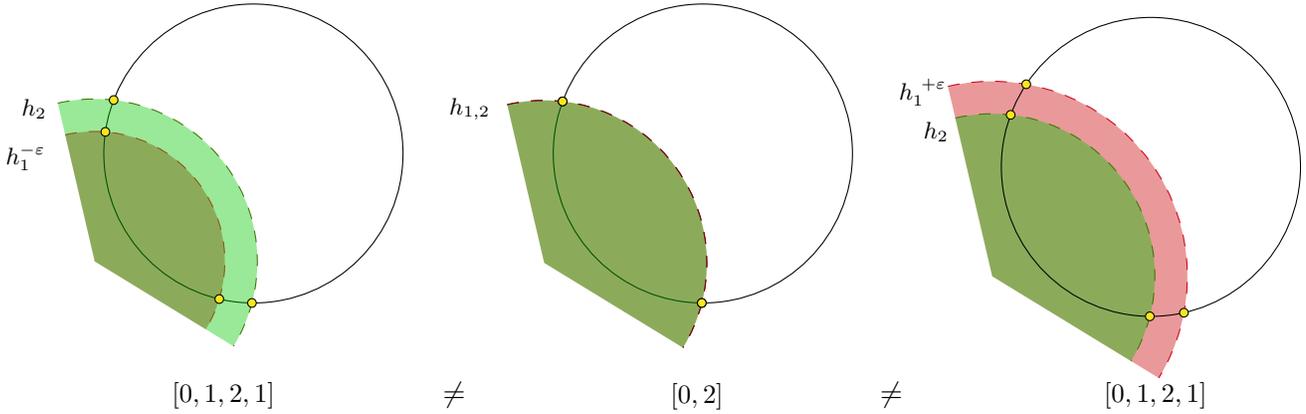

\end{Example}

We observe that in these prototypical examples, the cause of  non-robustness is either multiplicity of some boundary point ($\bdry h_i \cap \bdry h_j\neq \emptyset,\;i\neq j$) or the usage of a non-robust signal ($\card\bdry h_i\notin \{0, 2\}$) in (\ref{coupledsignal}). The next result states that this accounts for all possibilities, yielding a necessary and sufficient condition for robustness.

\subsection{Characterization of robustness}

\begin{Proposition}\label{p-robust} $H_n(v)$ is robust if and only if $\forall i,j$, $1\leq i<j\leq n$ we have  $\bdry h_i \cap \bdry h_j=\emptyset$ and $\card \bdry h_i\in \{0,2\}$.
\end{Proposition}

\textbf{Proof.} By Lemma~\ref{l-bdry-points}, the boundary points of $\bdry H_n(v)$ are exactly the discontinuity points of $H_n(v)$; the only exceptional scenario in Lemma~\ref{l-bdry-points} is when $h_i^{-1}(1)\supseteq S^1$, which is easily seen to fit into the non-robust case~(ii) of Lemma~\ref{l-element-counts-for-robustness} and hence is not relevant here. Definition~\ref{d-robust} of robustness asks that the sequence representations stay the same under sufficiently small perturbations: ${\text{seq}(H_n(v+\eps u))=\text{seq}(H_n(v))}$. With the equivalent sequences being defined by~\eqref{stringclass} and~\eqref{seqmapdefn}, it is clear that this is the case if and only if under such perturbations, boundary points in $\bdry H_n(v)$ do not disappear or merge or cross each other and new boundary points do not appear. The proof in both directions will rely on this observation.

\noindent (if) Suppose that the conditions in the statement are satisfied. By Lemma~\ref{l-element-counts-for-robustness}\,(i), for each $h_i$ we can find an $\bar\eps_i>0$ such that for all ${\bar\eps_i>\eps_i>0}, \lvert u\rvert=1$ we have $\card\bdry h_i(v+\eps_i u)=\card\bdry h_i (v)$. Additionally, it is easy to see that these small perturbations move the elements of each $\bdry h_i$ continuously. By assumption, each of the elements of $\bdry H_n(v)$ belongs to the boundary of a unique signal function $h_i$; let these elements be $x_j$, $j=1,\dots,\card \bdry H_n(v)$. We can find their neighborhoods $U_{x_j}\ni x_j$ such that
$U_{x_j} \cap U_{x_\ell}= \emptyset$ for $j\ne \ell$. By continuity, there exists an $\bar\eps\in(0,\min_{1\le i\le n}\bar\eps_i]$ such that for all ${\bar\eps>\eps>0}\,,\, {\lvert u \rvert=1}$, the perturbed decomposition $H_n(v+\eps u)$ has the property that the image of each original boundary point $x_j$ under the perturbation still belongs to $U_{x_j}$. We see that under such $\eps$-perturbations the number and order of the boundary points in $\bdry H_n(v)$ is preserved, and robustness follows.

\noindent (only if) We will use contraposition. Suppose that for some ${i\neq j}$ we have {${\bdry h_i \cap \bdry h_j \neq\emptyset}$ or ${\card \bdry h_i\notin \{0,2\}}$}. We want to show the decomposition to be not robust. If ${\card \bdry h_i \in \{1,\infty \}}$ for some~$i$, then there exist arbitrarily small perturbations of $h_i$ that cause vanishing/creation of a boundary point (Lemma~\ref{l-element-counts-for-robustness}\,(ii)). This changes the length of the sequence representation, meaning that $H_n(v)$ is non-robust.
Now suppose that for some ${i\neq j}$ we have ${\bdry h_i \cap\bdry h_j \neq\emptyset}$. We can assume that $\card \bdry h_i=\card\bdry h_j=2$ because otherwise we are done by the previous argument. It is easy to see that when $\card\bdry h_i=2$, for any neighborhood of $(m_i,r_i)$ we can always find a perturbation that changes the location of one boundary point only, keeping the other unchanged. Using such a perturbation, we can increase the number of boundary points of the decomposition, changing its sequence representation. Therefore, $H_n(v)$ is non-robust. \qed

The next corollary will be important for subsequent developments.

\begin{Corollary}\label{c-max-values}
If $H_n(v)$ is robust, then the maximum of the values of $H_n(v)$ around each $x\in S^1$ must equal the value at $x$, and the values around any $x\in S^1$ cannot have a difference larger than~1. More formally, if $H_{n}(v)$ is robust, then
\begin{eqnarray}
    & H_{n}(v)(x)=\max(H_{n}(v)(x^+),H_{n}(v)(x^-)\quad  \forall x \in S^1  \label{properUSC}\,,\\ & \max_{x\in S^1}\lvert H_{n}(v)(x^+)-H_{n}(v)(x^-)\rvert\leq 1\label{robuststep}
\end{eqnarray}
where $H_{n}(v)(x)$ is the evaluation of the function $H_{n}(v)$ at $x$ and $x^+,x^-$ are one-sided limits when approaching $x$ along the circle
 under a CCW parametrization.
\end{Corollary}

\textbf{Proof.} Assume that $H_n(v)$ is robust. We know from Lemma~\ref{l-bdry-points} that all discontinuity points of $H_n(v)$ are boundary points. For $x\in S^1 \setminus \bdry H_n(v)$, we must then have some neighborhood of $x$ for which $H_n(v)$ is constant, and~\eqref{properUSC} and~\eqref{robuststep} trivially hold. Now suppose that $x\in \bdry H_n(v)$. By Proposition~\ref{p-robust}, we must have a unique $h_i$ such that $x\in \bdry h_i$ with $\card \bdry h_i=2$. The intersection of the support of this $h_i$ with $S^1$ corresponds to an arc of the circle, thus there exists a left or right neighborhood of $x$ on which $h_i$ is identically 1. As $x\not\in \bdry h_j,\;j\neq i$, a small enough neighborhood of $x$ either entirely lies in the support of $h_j$ or is disjoint from it. Thus, all $h_{j},\; {j\neq i} ,$ are constant when restricted to some neighborhood of $x$, the only difference of 1 between $H_{n}(v)(x^+)$ and $H_{n}(v)(x^-)$ coming from $h_i$. Recalling that $h_i$ has a compact support by~\eqref{e-h-binary-def}, we conclude that $H_n(v)$ must satisfy~(\ref{properUSC}) and~(\ref{robuststep}). \qed

\section{Topological properties of decompositions}\label{s-top-prop}

Throughout this section, the following notion will be important.

\begin{Definition} A decomposition $H_n(v)$ is said to be~\textit{properly upper semicontinious} (properly u.s.c.) if it satisfies~(\ref{properUSC}).
\end{Definition}
Per Corollary~\ref{c-max-values}, robust $H_n(v)$ are properly u.s.c. Given an equivalence class of functions $\bar f $, we are now interested in the structure of all properly u.s.c.\:${H_n(v)\in \bar f}$. We must decide which types of signals $h_i$ (with $\card\bdry h_i=0,1,2$ or $\infty$) \textit{can} be used and if so, how many \textit{must} be used. For example: when decomposing functions taking strictly positive values everywhere on $S^1$, we could have some $h_i$ with $h_i^{-1}\supseteq S^1$ within (\ref{coupledsignal}). We can then ask, what is the smallest number of such $h_i$ we must use for $H_n(v)$ to be properly u.s.c.? (This question will be answered later in Proposition~\ref{p-baseline-levels}).

It is easy to see that every equivalence class $\bar f$ contains a unique properly u.s.c.\ function, which we denote $\langle \bar f\rangle$. This is the function whose only discontinuities are non-removable and whose values at these discontinuities are consistent with~\eqref{properUSC}. Using this ``canonical" representative of $\bar f$, we define the following quantities for the entire class $\bar f$:
\begin{eqnarray}
    &\bar f^*:= \max_{S^1}\langle \bar f\rangle,\qquad \bar f_* := \min_{S^1}\langle \bar f\rangle,\label{maxmindefn}\\ &
    \tau(\bar f) := \sum_{k=\bar f_*+1}^{\bar f^*}\big(\text{\#\{\small{disjoint sets making up} } \langle \bar f\rangle^{-1}([k,\infty))\} \,-1\big)\,.\label{taudefn}
\end{eqnarray}

We note that $\bar f^*$ could be equivalently defined as
$\min_{g\in \bar f}\,\max_{S^1}(g)$, with the minimum achieved at $\langle \bar f\rangle$, and similarly for $\bar f_*$. We will sometimes write simply $\tau$, omitting its argument $\bar f$ if it is clear from the context.
Letting
$$
\tau_k(\bar f):=\text{\#\{\small{disjoint sets making up} } \langle \bar f\rangle^{-1}([k,\infty))\}\,-1,\qquad k=0,\dots, \bar f^*
$$
we have $\tau=\sum_{k=\bar f_*+1}^{\bar f^*}\tau_k$.
Observe that this summation could equivalently be started from any integer $k'<\bar f_*+1$ because for $k\le \bar f_*$ we have $\langle \bar f\rangle^{-1}([k,\infty))=S^1$ hence $\tau_k=0$.

\begin{Example}\label{x-three-sets}
Consider the function $f$ (or, more precisely, its equivalence class $\bar f$) illustrated in Figure~\ref{f-example-3}. We want to find all properly u.s.c.~${H_n(v)\in\bar f}$. One can be constructed rather easily from the sets shown in the figure. To construct a complete set, however, we need to better understand what it means for $H_n(v)$ to be properly u.s.c. For example, is there a properly u.s.c.\ (in this case, actually {robust} as (\ref{robuststep}) is also satisfied) ${H_n(v)\in \bar f}$ that has no $h_i^{-1}(1)\supseteq S^1$ in (\ref{coupledsignal})? We will be able to resolve this question soon (see the end of Section~\ref{ss-char-rusc}).
\end{Example}

\begin{figure}[h!]
\tikzset{every picture/.style={line width=0.75pt}} 
\begin{tikzpicture}[x=0.75pt,y=0.75pt,yscale=-1,xscale=1]

\draw  [draw opacity=0] (442.79,199.07) .. controls (454.35,184.67) and (480.51,173) .. (501.24,173) .. controls (521.96,173) and (529.39,184.67) .. (517.83,199.07) .. controls (506.27,213.47) and (480.11,225.15) .. (459.38,225.15) .. controls (438.66,225.15) and (431.23,213.47) .. (442.79,199.07) -- cycle ;
\draw  [fill={rgb, 255:red, 65; green, 117; blue, 5 }  ,fill opacity=0.3 ] (383.11,186.09) -- (516.01,186.09) -- (441.89,278.45) -- (308.99,278.45) -- cycle ;
\draw  [color={rgb, 255:red, 65; green, 117; blue, 5 }  ,draw opacity=1 ][line width=2.25]  (359.61,232.27) .. controls (375.9,211.97) and (412.79,195.51) .. (442,195.51) .. controls (471.21,195.51) and (481.69,211.97) .. (465.4,232.27) .. controls (449.1,252.57) and (412.21,269.03) .. (383,269.03) .. controls (353.79,269.03) and (343.31,252.57) .. (359.61,232.27) -- cycle ;
\draw  [draw opacity=0][line width=0.75]  (333.37,216.76) .. controls (343.81,203.76) and (367.44,193.22) .. (386.15,193.22) .. controls (404.87,193.22) and (411.58,203.76) .. (401.14,216.76) .. controls (390.7,229.77) and (367.07,240.31) .. (348.36,240.31) .. controls (329.64,240.31) and (322.93,229.77) .. (333.37,216.76) -- cycle ;
\draw  [draw opacity=0] (356.12,266.27) .. controls (367.52,252.06) and (393.34,240.54) .. (413.79,240.54) .. controls (434.24,240.54) and (441.57,252.06) .. (430.17,266.27) .. controls (418.76,280.48) and (392.94,292) .. (372.49,292) .. controls (352.04,292) and (344.71,280.48) .. (356.12,266.27) -- cycle ;
\draw  [draw opacity=0] (384.82,263.05) .. controls (390.01,256.58) and (401.77,251.34) .. (411.08,251.34) .. controls (420.39,251.34) and (423.73,256.58) .. (418.54,263.05) .. controls (413.34,269.52) and (401.58,274.77) .. (392.27,274.77) .. controls (382.96,274.77) and (379.62,269.52) .. (384.82,263.05) -- cycle ;
\draw  [draw opacity=0] (455.57,203.99) .. controls (460.1,198.34) and (470.36,193.77) .. (478.49,193.77) .. controls (486.61,193.77) and (489.53,198.34) .. (484.99,203.99) .. controls (480.46,209.64) and (470.2,214.21) .. (462.08,214.21) .. controls (453.95,214.21) and (451.04,209.64) .. (455.57,203.99) -- cycle ;
\draw  [draw opacity=0] (380.82,180.05) .. controls (386.01,173.58) and (397.77,168.34) .. (407.08,168.34) .. controls (416.39,168.34) and (419.73,173.58) .. (414.54,180.05) .. controls (409.34,186.52) and (397.58,191.77) .. (388.27,191.77) .. controls (378.96,191.77) and (375.62,186.52) .. (380.82,180.05) -- cycle ;
\draw  [color={rgb, 255:red, 0; green, 0; blue, 0 }  ,draw opacity=1 ][fill={rgb, 255:red, 208; green, 2; blue, 27 }  ,fill opacity=0.3 ] (379.11,103.09) -- (512.01,103.09) -- (437.89,195.45) -- (304.99,195.45) -- cycle ;
\draw  [dash pattern={on 4.5pt off 4.5pt}][line width=0.75]  (355.61,149.27) .. controls (371.9,128.97) and (408.79,112.51) .. (438,112.51) .. controls (467.21,112.51) and (477.69,128.97) .. (461.4,149.27) .. controls (445.1,169.57) and (408.21,186.03) .. (379,186.03) .. controls (349.79,186.03) and (339.31,169.57) .. (355.61,149.27) -- cycle ;
\draw  [draw opacity=0][line width=0.75]  (329.37,133.76) .. controls (339.81,120.76) and (363.44,110.22) .. (382.15,110.22) .. controls (400.87,110.22) and (407.58,120.76) .. (397.14,133.76) .. controls (386.7,146.77) and (363.07,157.31) .. (344.36,157.31) .. controls (325.64,157.31) and (318.93,146.77) .. (329.37,133.76) -- cycle ;
\draw  [draw opacity=0] (438.79,116.07) .. controls (450.35,101.67) and (476.51,90) .. (497.24,90) .. controls (517.96,90) and (525.39,101.67) .. (513.83,116.07) .. controls (502.27,130.47) and (476.11,142.15) .. (455.38,142.15) .. controls (434.66,142.15) and (427.23,130.47) .. (438.79,116.07) -- cycle ;
\draw  [draw opacity=0] (352.12,183.27) .. controls (363.52,169.06) and (389.34,157.54) .. (409.79,157.54) .. controls (430.24,157.54) and (437.57,169.06) .. (426.17,183.27) .. controls (414.76,197.48) and (388.94,209) .. (368.49,209) .. controls (348.04,209) and (340.71,197.48) .. (352.12,183.27) -- cycle ;
\draw  [draw opacity=0][line width=2.25]  (350.84,156.47) .. controls (352.05,154.14) and (353.63,151.73) .. (355.61,149.27) .. controls (365.17,137.35) and (381.83,126.76) .. (399.87,120.05) -- (408.5,149.27) -- cycle ; \draw  [color={rgb, 255:red, 208; green, 2; blue, 27 }  ,draw opacity=1 ][line width=2.25]  (350.84,156.47) .. controls (352.05,154.14) and (353.63,151.73) .. (355.61,149.27) .. controls (365.17,137.35) and (381.83,126.76) .. (399.87,120.05) ;
\draw  [draw opacity=0][line width=2.25]  (441.41,112.59) .. controls (465.38,113.71) and (475.08,126.77) .. (465.31,143.61) -- (408.5,149.27) -- cycle ; \draw  [color={rgb, 255:red, 208; green, 2; blue, 27 }  ,draw opacity=1 ][line width=2.25]  (441.41,112.59) .. controls (465.38,113.71) and (475.08,126.77) .. (465.31,143.61) ;
\draw  [draw opacity=0][line width=2.25]  (430.88,172.48) .. controls (414.55,180.75) and (395.54,186.03) .. (379,186.03) .. controls (368.81,186.03) and (360.89,184.02) .. (355.6,180.55) -- (408.5,149.27) -- cycle ; \draw  [color={rgb, 255:red, 208; green, 2; blue, 27 }  ,draw opacity=1 ][line width=2.25]  (430.88,172.48) .. controls (414.55,180.75) and (395.54,186.03) .. (379,186.03) .. controls (368.81,186.03) and (360.89,184.02) .. (355.6,180.55) ;
\draw  [dash pattern={on 4.5pt off 4.5pt}][line width=0.75]  (349.61,65.27) .. controls (365.9,44.97) and (402.79,28.51) .. (432,28.51) .. controls (461.21,28.51) and (471.69,44.97) .. (455.4,65.27) .. controls (439.1,85.57) and (402.21,102.03) .. (373,102.03) .. controls (343.79,102.03) and (333.31,85.57) .. (349.61,65.27) -- cycle ;
\draw  [line width=1.5]  (96,152) .. controls (96,107.82) and (131.82,72) .. (176,72) .. controls (220.18,72) and (256,107.82) .. (256,152) .. controls (256,196.18) and (220.18,232) .. (176,232) .. controls (131.82,232) and (96,196.18) .. (96,152) -- cycle ;
\draw  [color={rgb, 255:red, 255; green, 255; blue, 255 }  ,draw opacity=1 ][line width=0.75]  (37.5,118.25) .. controls (37.5,89.95) and (60.45,67) .. (88.75,67) .. controls (117.05,67) and (140,89.95) .. (140,118.25) .. controls (140,146.55) and (117.05,169.5) .. (88.75,169.5) .. controls (60.45,169.5) and (37.5,146.55) .. (37.5,118.25) -- cycle ;
\draw  [draw opacity=0] (181.5,79.75) .. controls (181.5,48.41) and (206.91,23) .. (238.25,23) .. controls (269.59,23) and (295,48.41) .. (295,79.75) .. controls (295,111.09) and (269.59,136.5) .. (238.25,136.5) .. controls (206.91,136.5) and (181.5,111.09) .. (181.5,79.75) -- cycle ;
\draw  [draw opacity=0] (132,226) .. controls (132,195.07) and (157.07,170) .. (188,170) .. controls (218.93,170) and (244,195.07) .. (244,226) .. controls (244,256.93) and (218.93,282) .. (188,282) .. controls (157.07,282) and (132,256.93) .. (132,226) -- cycle ;
\draw  [draw opacity=0] (171.5,219) .. controls (171.5,204.92) and (182.92,193.5) .. (197,193.5) .. controls (211.08,193.5) and (222.5,204.92) .. (222.5,219) .. controls (222.5,233.08) and (211.08,244.5) .. (197,244.5) .. controls (182.92,244.5) and (171.5,233.08) .. (171.5,219) -- cycle ;
\draw  [draw opacity=0] (209,93.75) .. controls (209,81.46) and (218.96,71.5) .. (231.25,71.5) .. controls (243.54,71.5) and (253.5,81.46) .. (253.5,93.75) .. controls (253.5,106.04) and (243.54,116) .. (231.25,116) .. controls (218.96,116) and (209,106.04) .. (209,93.75) -- cycle ;
\draw  [draw opacity=0] (374.82,96.05) .. controls (380.01,89.58) and (391.77,84.34) .. (401.08,84.34) .. controls (410.39,84.34) and (413.73,89.58) .. (408.54,96.05) .. controls (403.34,102.52) and (391.58,107.77) .. (382.27,107.77) .. controls (372.96,107.77) and (369.62,102.52) .. (374.82,96.05) -- cycle ;
\draw  [draw opacity=0][line width=2.25]  (408.88,95.31) .. controls (396.86,99.54) and (384.35,102.03) .. (373,102.03) .. controls (372.86,102.03) and (372.72,102.03) .. (372.58,102.03) -- (402.5,65.27) -- cycle ; \draw  [color={rgb, 255:red, 162; green, 2; blue, 193 }  ,draw opacity=1 ][line width=2.25]  (408.88,95.31) .. controls (396.86,99.54) and (384.35,102.03) .. (373,102.03) .. controls (372.86,102.03) and (372.72,102.03) .. (372.58,102.03) ;
\draw  [draw opacity=0] (451.57,120.99) .. controls (456.1,115.34) and (466.36,110.77) .. (474.49,110.77) .. controls (482.61,110.77) and (485.53,115.34) .. (480.99,120.99) .. controls (476.46,126.64) and (466.2,131.21) .. (458.08,131.21) .. controls (449.95,131.21) and (447.04,126.64) .. (451.57,120.99) -- cycle ;
\draw  [fill={rgb, 255:red, 144; green, 19; blue, 254 }  ,fill opacity=0.3 ] (373.11,19.09) -- (506.01,19.09) -- (431.89,111.45) -- (298.99,111.45) -- cycle ;
\draw  [draw opacity=0][line width=0.75]  (323.37,49.76) .. controls (333.81,36.76) and (357.44,26.22) .. (376.15,26.22) .. controls (394.87,26.22) and (401.58,36.76) .. (391.14,49.76) .. controls (380.7,62.77) and (357.07,73.31) .. (338.36,73.31) .. controls (319.64,73.31) and (312.93,62.77) .. (323.37,49.76) -- cycle ;
\draw  [draw opacity=0] (432.79,32.07) .. controls (444.35,17.67) and (470.51,6) .. (491.24,6) .. controls (511.96,6) and (519.39,17.67) .. (507.83,32.07) .. controls (496.27,46.47) and (470.11,58.15) .. (449.38,58.15) .. controls (428.66,58.15) and (421.23,46.47) .. (432.79,32.07) -- cycle ;
\draw  [draw opacity=0] (346.12,99.27) .. controls (357.52,85.06) and (383.34,73.54) .. (403.79,73.54) .. controls (424.24,73.54) and (431.57,85.06) .. (420.17,99.27) .. controls (408.76,113.48) and (382.94,125) .. (362.49,125) .. controls (342.04,125) and (334.71,113.48) .. (346.12,99.27) -- cycle ;
\draw  [draw opacity=0] (445.57,36.99) .. controls (450.1,31.34) and (460.36,26.77) .. (468.49,26.77) .. controls (476.61,26.77) and (479.53,31.34) .. (474.99,36.99) .. controls (470.46,42.64) and (460.2,47.21) .. (452.08,47.21) .. controls (443.95,47.21) and (441.04,42.64) .. (445.57,36.99) -- cycle ;
\draw  [draw opacity=0][line width=2.25]  (451.59,31.94) .. controls (458.22,34.85) and (462.12,39.45) .. (462.92,45.16) -- (402.5,65.27) -- cycle ; \draw  [color={rgb, 255:red, 162; green, 2; blue, 193 }  ,draw opacity=1 ][line width=2.25]  (451.59,31.94) .. controls (458.22,34.85) and (462.12,39.45) .. (462.92,45.16) ;

\draw (103.6,124.33) node [anchor=north west][inner sep=0.75pt]  [font=\small,color={rgb, 255:red, 208; green, 2; blue, 27 }  ,opacity=1 ,rotate=-292.65]  {$4$};
\draw (148.52,80.47) node [anchor=north west][inner sep=0.75pt]  [font=\small,color={rgb, 255:red, 65; green, 117; blue, 5 }  ,opacity=1 ,rotate=-345.32]  {$3$};
\draw (121.09,181.81) node [anchor=north west][inner sep=0.75pt]  [font=\small,color={rgb, 255:red, 65; green, 117; blue, 5 }  ,opacity=1 ,rotate=-49.28]  {$3$};
\draw (234.5,174.25) node [anchor=north west][inner sep=0.75pt]  [font=\small,color={rgb, 255:red, 65; green, 117; blue, 5 }  ,opacity=1 ,rotate=-288.61]  {$3$};
\draw (153.22,210.59) node [anchor=north west][inner sep=0.75pt]  [font=\small,color={rgb, 255:red, 208; green, 2; blue, 27 }  ,opacity=1 ,rotate=-17.18]  {$4$};
\draw (242.76,119.8) node [anchor=north west][inner sep=0.75pt]  [font=\small,color={rgb, 255:red, 208; green, 2; blue, 27 }  ,opacity=1 ,rotate=-65.47]  {$4$};
\draw (191.24,77.08) node [anchor=north west][inner sep=0.75pt]  [font=\small,color={rgb, 255:red, 208; green, 2; blue, 27 }  ,opacity=1 ,rotate=-14.94]  {$4$};
\draw (224.32,95.5) node [anchor=north west][inner sep=0.75pt]  [font=\small,color={rgb, 255:red, 162; green, 2; blue, 193 }  ,opacity=1 ,rotate=-40.06]  {$5$};
\draw (186.98,212.7) node [anchor=north west][inner sep=0.75pt]  [font=\small,color={rgb, 255:red, 162; green, 2; blue, 193 }  ,opacity=1 ,rotate=-346.72]  {$5$};
\draw (214.79,199.16) node [anchor=north west][inner sep=0.75pt]  [font=\small,color={rgb, 255:red, 208; green, 2; blue, 27 }  ,opacity=1 ,rotate=-317.7]  {$4$};
\draw (514.01,106.49) node [anchor=north west][inner sep=0.75pt]    {$^{f^{-1}([ 4,\infty ))}$};
\draw (508.01,22.49) node [anchor=north west][inner sep=0.75pt]    {$^{f^{-1}([ 5,\infty ))}$};
\draw (518.01,189.49) node [anchor=north west][inner sep=0.75pt]    {$^{f^{-1}([ k,\infty)) \ \ 1\leq k\leq 3\ }$};
\draw (454.58,264.51) node [anchor=north west][inner sep=0.75pt]  [font=\footnotesize]  {$\tau _{k} =0$};
\draw (133,49.4) node [anchor=north west][inner sep=0.75pt]  [font=\small]  {$f:S^{1}\rightarrow \mathbb{N}$};
\draw (306.99,198.85) node [anchor=north west][inner sep=0.75pt]  [font=\footnotesize]  {$\tau _{4} =2$};
\draw (300.99,114.85) node [anchor=north west][inner sep=0.75pt]  [font=\footnotesize]  {$\tau _{5} =1$};

\draw [fill={rgb, 255:red, 248; green, 231; blue, 28 }  ,fill opacity=1 ]  (129.37, 86.99) circle [x radius= 2, y radius= 2]   ;
\draw [fill={rgb, 255:red, 248; green, 231; blue, 28 }  ,fill opacity=1 ]  (97.75, 168.71) circle [x radius= 2, y radius= 2]   ;
\draw [fill={rgb, 255:red, 248; green, 231; blue, 28 }  ,fill opacity=1 ]  (181.99, 72.22) circle [x radius= 2, y radius= 2]   ;
\draw [fill={rgb, 255:red, 248; green, 231; blue, 28 }  ,fill opacity=1 ]  (254.03, 134.28) circle [x radius= 2, y radius= 2]   ;
\draw [fill={rgb, 255:red, 248; green, 231; blue, 28 }  ,fill opacity=1 ]  (132.42, 219.1) circle [x radius= 2, y radius= 2]   ;
\draw [fill={rgb, 255:red, 248; green, 231; blue, 28 }  ,fill opacity=1 ]  (238.55, 201.88) circle [x radius= 2, y radius= 2]   ;
\draw [fill={rgb, 255:red, 248; green, 231; blue, 28 }  ,fill opacity=1 ]  (222.43, 217.15) circle [x radius= 2, y radius= 2]   ;
\draw [fill={rgb, 255:red, 248; green, 231; blue, 28 }  ,fill opacity=1 ]  (175.05, 231.99) circle [x radius= 2, y radius= 2]   ;
\draw [fill={rgb, 255:red, 248; green, 231; blue, 28 }  ,fill opacity=1 ]  (212.98, 81.04) circle [x radius= 2, y radius= 2]   ;
\draw [fill={rgb, 255:red, 248; green, 231; blue, 28 }  ,fill opacity=1 ]  (244.9, 111.32) circle [x radius= 2, y radius= 2]   ;
\draw [fill={rgb, 255:red, 248; green, 231; blue, 28 }  ,fill opacity=1 ]  (419.17, 262.21) circle [x radius= 2, y radius= 2]   ;
\draw [fill={rgb, 255:red, 248; green, 231; blue, 28 }  ,fill opacity=1 ]  (382.37, 269.02) circle [x radius= 2, y radius= 2]   ;
\draw [fill={rgb, 255:red, 248; green, 231; blue, 28 }  ,fill opacity=1 ]  (445.89, 195.62) circle [x radius= 2, y radius= 2]   ;
\draw [fill={rgb, 255:red, 248; green, 231; blue, 28 }  ,fill opacity=1 ]  (470.63, 224.13) circle [x radius= 2, y radius= 2]   ;
\draw [fill={rgb, 255:red, 248; green, 231; blue, 28 }  ,fill opacity=1 ]  (405.64, 202.4) circle [x radius= 2, y radius= 2]   ;
\draw [fill={rgb, 255:red, 248; green, 231; blue, 28 }  ,fill opacity=1 ]  (354.6, 239.95) circle [x radius= 2, y radius= 2]   ;
\draw [fill={rgb, 255:red, 248; green, 231; blue, 28 }  ,fill opacity=1 ]  (435.46, 255.19) circle [x radius= 2, y radius= 2]   ;
\draw [fill={rgb, 255:red, 248; green, 231; blue, 28 }  ,fill opacity=1 ]  (358.94, 263.1) circle [x radius= 2, y radius= 2]   ;
\draw [fill={rgb, 255:red, 248; green, 231; blue, 28 }  ,fill opacity=1 ]  (461.56, 198.93) circle [x radius= 2, y radius= 2]   ;
\draw [fill={rgb, 255:red, 248; green, 231; blue, 28 }  ,fill opacity=1 ]  (472.91, 212.04) circle [x radius= 2, y radius= 2]   ;
\draw [fill={rgb, 255:red, 248; green, 231; blue, 28 }  ,fill opacity=1 ]  (415.17, 179.21) circle [x radius= 2, y radius= 2]   ;
\draw [fill={rgb, 255:red, 248; green, 231; blue, 28 }  ,fill opacity=1 ]  (378.37, 186.02) circle [x radius= 2, y radius= 2]   ;
\draw [fill={rgb, 255:red, 248; green, 231; blue, 28 }  ,fill opacity=1 ]  (401.64, 119.4) circle [x radius= 2, y radius= 2]   ;
\draw [fill={rgb, 255:red, 248; green, 231; blue, 28 }  ,fill opacity=1 ]  (350.6, 156.95) circle [x radius= 2, y radius= 2]   ;
\draw [fill={rgb, 255:red, 248; green, 231; blue, 28 }  ,fill opacity=1 ]  (441.89, 112.62) circle [x radius= 2, y radius= 2]   ;
\draw [fill={rgb, 255:red, 248; green, 231; blue, 28 }  ,fill opacity=1 ]  (466.63, 141.13) circle [x radius= 2, y radius= 2]   ;
\draw [fill={rgb, 255:red, 248; green, 231; blue, 28 }  ,fill opacity=1 ]  (431.46, 172.19) circle [x radius= 2, y radius= 2]   ;
\draw [fill={rgb, 255:red, 248; green, 231; blue, 28 }  ,fill opacity=1 ]  (354.94, 180.1) circle [x radius= 2, y radius= 2]   ;
\draw [fill={rgb, 255:red, 248; green, 231; blue, 28 }  ,fill opacity=1 ]  (457.56, 115.93) circle [x radius= 2, y radius= 2]   ;
\draw [fill={rgb, 255:red, 248; green, 231; blue, 28 }  ,fill opacity=1 ]  (468.91, 129.04) circle [x radius= 2, y radius= 2]   ;
\draw [fill={rgb, 255:red, 248; green, 231; blue, 28 }  ,fill opacity=1 ]  (409.17, 95.21) circle [x radius= 2, y radius= 2]   ;
\draw [fill={rgb, 255:red, 248; green, 231; blue, 28 }  ,fill opacity=1 ]  (372.37, 102.02) circle [x radius= 2, y radius= 2]   ;
\draw [fill={rgb, 255:red, 248; green, 231; blue, 28 }  ,fill opacity=1 ]  (395.64, 35.4) circle [x radius= 2, y radius= 2]   ;
\draw [fill={rgb, 255:red, 248; green, 231; blue, 28 }  ,fill opacity=1 ]  (344.6, 72.95) circle [x radius= 2, y radius= 2]   ;
\draw [fill={rgb, 255:red, 248; green, 231; blue, 28 }  ,fill opacity=1 ]  (435.89, 28.62) circle [x radius= 2, y radius= 2]   ;
\draw [fill={rgb, 255:red, 248; green, 231; blue, 28 }  ,fill opacity=1 ]  (460.63, 57.13) circle [x radius= 2, y radius= 2]   ;
\draw [fill={rgb, 255:red, 248; green, 231; blue, 28 }  ,fill opacity=1 ]  (425.46, 88.19) circle [x radius= 2, y radius= 2]   ;
\draw [fill={rgb, 255:red, 248; green, 231; blue, 28 }  ,fill opacity=1 ]  (348.94, 96.1) circle [x radius= 2, y radius= 2]   ;
\draw [fill={rgb, 255:red, 248; green, 231; blue, 28 }  ,fill opacity=1 ]  (451.56, 31.93) circle [x radius= 2, y radius= 2]   ;
\draw [fill={rgb, 255:red, 248; green, 231; blue, 28 }  ,fill opacity=1 ]  (462.91, 45.04) circle [x radius= 2, y radius= 2]   ;
\end{tikzpicture}
\caption{${f\in \bar f}$ in Example~\ref{x-three-sets} and a properly u.s.c $H_n(v)$ using three signals with $h_i^{-1}(1)\supseteq S^1$.\label{f-example-3}}
\end{figure}

Given $H_n(v)\in \bar f$, we define the chain of \textit{upper excursion sets} \cite{euler}:
\begin{equation}\label{upperexcursionchain}
    H_n(v)^{-1}([\bar f^*,\infty))\subseteq \dots \subseteq H_n(v)^{-1}([\bar f_*+2,\infty)) \subseteq  H_n(v)^{-1}([\bar f_*+1,\infty))\,.
\end{equation}

\begin{Proposition}\label{p-discont} If ${H_n(v)\in \bar f}$ is robust, then $\card\bdry H_n(v)=2(\bar f^*-\bar f_*+\tau)$.
\end{Proposition}

\textbf{Proof.} In view of Lemma~\ref{l-bdry-points} and Proposition~\ref{p-robust}, we see that {every} $x\in \bdry H_n(v)$ is a discontinuity point of the robust ${H_n(v)\in\bar f}$. Therefore, it suffices to count the discontinuity points of $H_n(v)$. Since there is a unique function $\langle\bar f\rangle\in \bar f$ that is properly u.s.c., any two robust (hence properly u.s.c.) decompositions in $\bar f$ must be equal as {functions}, and so the sets in~\eqref{upperexcursionchain} are uniquely defined.
By~(\ref{robuststep}), the boundary points of the different sets appearing in~(\ref{upperexcursionchain}) must be disjoint. Computing $\card\bdry H_n(v)$ thus reduces to counting the boundary points of each element of (\ref{upperexcursionchain}) and summing. Furthermore, every element of the chain (\ref{upperexcursionchain}) must be made up of a union of disjoint sets with two boundary points each; indeed, otherwise it will contain an isolated point and we will have a contradiction to (\ref{properUSC}). Therefore, $\card\bdry H_n(v)$ must equal twice the number of disjoint sets making up the chain (\ref{upperexcursionchain}), and this is readily seen to be $2(\tau + \bar f^*- \bar f_*)$ as claimed. \qed

\subsection{Characterization of properly u.s.c.\ decompositions}\label{ss-char-rusc}

\begin{Lemma}\label{l-pusc-char} $H_n(v)\in \bar f$ is properly u.s.c.\ if and only if $\:\forall i,j$ such that $1\leq i\leq j\leq n$,
if $\bdry h_i \cap  \bdry h_j\neq \emptyset$, then no point $x\in \bdry h_i\cap \bdry h_j$ is an isolated point of $h_i^{-1}(1)\cap h_j^{-1}(1)\cap S^1$.
\end{Lemma}

The condition in the lemma is equivalent to: If  $x\in \bdry h_i\cap \bdry h_j$ with $\card \bdry h_i=\card\bdry h_j=2$ and $i\ne j$, then $h_i^{-1}(1)\cap S^1$ and $h_j^{-1}(1)\cap S^1$ extend in the same direction from $x$, and if $\card \bdry h_i = 1$ then $h_i^{-1}(1)\supseteq S^1$.

\textbf{Proof.} Since all $h_i$ are compactly supported in $S^1$, the corresponding sum (\ref{coupledsignal}) is {upper semicontinuous}\footnote{In the usual sense that the value of the function at any point is {larger} or equal to its maximal value in a sufficiently small neighborhood of the point, the equality in (\ref{properUSC}) being replaced by $\geq$.}, and properly u.s.c.\ away from the discontinuity points of the decomposition. Therefore, we only need to check discontinuity points of $H_n(v)$ to show that it is properly u.s.c.

\noindent (if) By Lemma~\ref{l-bdry-points}, the only discontinuity points of $H_n(v)$ are the
boundary points of $h_i$ with $h_i^{-1}(1)\not\supseteq S^1$. This condition rules out signals $h_i$ with $\card \bdry h_i=\infty$ and, by assumption, also those with $\card \bdry h_i=1$. Thus we need to only check if the value of $H_n(v)$ at $x\in\bdry h_i$ with $\card\bdry h_i=2$ is the maximum of the values around $x$. Then the application of our assumption to all $h_{i_j}$ with $x\in \bdry h_{i_j}$ shows that $x$ is not an isolated point of $h_{i_1}^{-1}(1)\cap\dots\cap h_{i_k}^{-1}(1)\cap S^1$, meaning that the point $x$ is an endpoint of some arc of the circle on which all $h_{i_j}$ are identically 1. For each of the remaining signals $h_j$, with $x\not \in \bdry h_j$, we always have a neighborhood of $x$ where $h_j$ is constant. Therefore, for each $x\in \bdry h_i$ with $\card\bdry h_i=2$, the value of $H_n(v)$ at $x$ is equal to the value of $H_n(v)$ on some one-sided neighborhood of $x$, which is the maximum in~\eqref{properUSC}.

\noindent(only if) We will use contraposition. Assume that $\exists i,j$ such that $x\in \bdry h_i\cap \bdry h_j$ is an isolated point of
$h_i^{-1}(1)\cap h_j^{-1}(1)\cap S^1$. If $i=j$, then we must have $\card\bdry  h_i=1$ and $h_i^{-1}(1)\not \supseteq S^1$. Such an $h_i$ takes value 1 only at a single point $x\in \bdry h_i$, and is 0 everywhere else in $S^1$. Because every $h_i$ in (\ref{coupledsignal}) is compactly supported, $H_n(v)$ cannot satisfy (\ref{properUSC}) at $x$, as its value is going to be larger (by 1) than the maximum value in a sufficiently small neighborhood of $x$. Now assume $i\neq j$; then we must have $\card\bdry h_i=\card\bdry h_j=2$ with supports of $h_i$ and $h_j$ extending in opposite directions around $x$ within $S^1$. For the point of intersection we can again argue, as in the previous case, that the values of $H_n(v)$ on any sufficiently small neighborhood of $x\in \bdry h_i\cap \bdry h_j$ will be smaller (by 1) than the value at $x$, violating~\eqref{properUSC}.~\qed

We infer from Lemma~\ref{l-pusc-char} (and its proof) that
the only discontinuity points of a properly u.s.c.\ $H_n(v)$ are the boundary points of $h_i$ with $\card \bdry h_i=2$. We also note the following simple fact.

\begin{Lemma}\label{l-tauinvariance}
    If $H_n(v)\in \bar f$ is properly u.s.c.\ with $\card\bdry h_i\in \{0,2\}\,\forall i$, then $
\exists \bar\eps>0\; \text{\,s.t.}\;\forall \eps \in (0,\bar\eps) \;\forall u\in \mathbb{R}^{3n} \;\text{with}\;\lvert u\rvert=1,\; \tau_k (\overline{H_{n}(v+\eps u)}) = \tau_k (\overline{H_{n}(v)})\; \text{for each } k=0,\dots, \bar f^*.$~$\omitted$
\end{Lemma}

The next proposition elaborates upon the above observations.

\begin{Proposition}\label{p-pusc-char} If $H_n(v)\in\bar f$ is properly u.s.c., then every $x\in S^1$ belongs to $\bdry h_i$ for exactly $k_x$ signal functions $h_i$ with $\card\bdry h_i=2$, where
$
k_x:= \lvert H_{n}(v)(x^+)-H_{n}(v)(x^-)\rvert.
$
Also,  $\sum_{x\in \bdry H_n(v)} k_x = 2(\bar f^*-\bar f_*+ \tau)$.
\end{Proposition}

\textbf{Proof. } Assume that $H_n(v)\in \bar f$ is properly u.s.c. Then, as stated before the proposition, any discontinuity point $x\in \bdry H_n(v)$ must be an element of some $\bdry h_i$ with $\card \bdry h_i=2$ only, and it is easy to see that there must be at least $k_x$ of these $\bdry h_i$. If $x$ is an element of more than $k_x$ many $\bdry h_i$ with $\card\bdry h_i=2$, however, we must have some $h_i,h_j$ that have supports locally intersecting only at $x$, implying that $x$ is an isolated point of $h_i^{-1}(1)\cap h_j^{-1}(1)\cap S^1$, a contradiction to Lemma~\ref{l-pusc-char}.

Now consider the chain of sets (\ref{upperexcursionchain}) again. Some boundary points of the sets ${H_n(v)}^{-1}([k,\infty))$ may not be disjoint because the sets must overlap, by definition, $k_x$ times at $x\in S^1$. Any $h_i$ with $\card \bdry h_i\not\in \{0,2\}$ in (\ref{coupledsignal}) must have $h_i^{-1}(1)\supseteq S^1$ by Lemma~\ref{l-pusc-char}, meaning that all such $h_i$ are identically 1 on $S^1$. We can then subtract such $h_i$ from (\ref{coupledsignal}) and this will not change $k_x$ for any $x\in S^1$. After this subtraction, we obtain a new properly u.s.c.\ decomposition $H_{n_1}(v_1)$, $n_1 \leq n$ which satisfies $\sum_{x\in \bdry H_n(v)}k_x=\sum_{x\in \bdry H_{n_1}({v_1})}k_x$ and only has $\card \bdry h_i\in \{0,2\}$. It is easy to show that ${\overline{H_{n_1}(v_1)}\,}^*=\bar f^*-(n-n_1)$, ${\overline{H_{n_1}(v_1)}\,}_*=\bar f_*-(n-n_1)$, and $\tau(\overline{H_{n_1}(v_1)})=\tau(\bar f)$.
Now, small perturbations to the $h_i$ of $H_{n_1}(v_1)$ cannot cause creation/vanishing of boundary points (Lemma~\ref{l-element-counts-for-robustness}\,(i)), hence they cannot change $\sum_{x\in \bdry H_{n_1}({v_1})} k_x$. Consider some perturbed decomposition $H_{n_1}(v_1+\eps u)$ that has no shared boundary points, i.e., $k_x=1$ for all $x\in \bdry H_{n_1}(v_1+\eps u)$. Such a perturbation can be made arbitrarily small, and by Proposition~\ref{p-robust} the new decomposition $H_{n_1}(v_1+\eps u)$ is robust.  Therefore, $\sum_{x\in \bdry H_{n_1}(v_1)} k_x$ must equal  $\sum_{x\in \bdry H_{n_1}(v_1+\eps u)} k_x$, which by Proposition~\ref{p-discont} is equal to  $2(\bar f^*-(n-n_1)-(\bar f_*-(n-n_1)+\tau(\overline{H_{n_1}(v_1+\eps u)})))=2(\bar f^*-\bar f_*+\tau(\overline{H_{n_1}(v_1+\eps u)}))=2(\bar f^*-\bar f_*+\tau(\overline{H_{n_1}(v_1)}))=2(\bar f^*-\bar f_*+\tau(\bar f))$, where the second equality relies on Lemma~\ref{l-tauinvariance}.~\qed

The following fact was established in the course of the previous proof and we state it as a separate corollary.

\begin{Corollary}\label{c-intervals} If ${H_n(v)\in \bar f}$ is properly u.s.c., then it must have exactly $\bar f^*-\bar f_*+\tau$ signal functions $h_i$ with $\card\bdry h_i=2$.
\end{Corollary}

\begin{Remark}\label{c-intervals_lowerbound}
An arbitrary (not necessarily properly u.s.c.) decomposition $H_n(v)\in \bar f$ has at least  $\bar f^*-\bar f_*+\tau$ signal functions $h_i$ with $\card \bdry h_i=2$. This can be seen by repeating the arguments of the previous proof, with each $x\in S^1$ necessarily belonging to \emph{at least} $k_x$ signal functions $h_i$ with $\card \bdry h_i=2$.
\end{Remark}

\begin{Proposition}\label{p-baseline-levels} If $H_n(v)\in\bar f$ is properly u.s.c., then it has at least $\max\{\bar f_*-\tau(\bar f),0\}$ signal functions $h_i$  with  $h_i^{-1}(1)\supseteq S^1$, i.e., $h_i\equiv 1$ on $S^1$.
\end{Proposition}

\textbf{Proof.} Let ${H_{n}(v)\in \bar f}$ be properly u.s.c. and suppose that it has $k<\bar f_*-\tau$ many $h_i$ with $h_i^{-1}(1)\supseteq S^1$. As $\bar f_*>k$, we have: $(\overline{f-k})^*=\bar f^*-k$, $(\overline{f-k})_*=\bar f_*-k$, and $\tau(\overline{f-k}) =\tau(\bar f)$. Therefore, by Corollary~\ref{c-intervals}, any properly u.s.c.\ ${H_{n_1}(v_1)\in \overline{f-k}}$ has $(\bar f^*-k)-(\bar f_*-k)+\tau=\bar f^*-\bar f_*+\tau$ many $h_i$, all with $\card\bdry h_i=2$ (we ignore any $h_i$ that are identically 0 on $S^1$).
The largest possible value of $H_{n_1}(v_1)$, attained at a point of intersection of
all $h_i$ (if one exists), is $\bar f^*-\bar f_*+\tau$. But we have $\bar f^*-\bar f_*+\tau<\bar f^*-k$, so this largest value is still smaller than the required maximum---a contradiction.~\qed

We now revisit Example~\ref{x-three-sets} and give another decomposition for the function considered there, this time one that uses no signals $h_i$ with $h_i^{-1}(1)\supseteq S^1$. Notice that $\bar f_*=\tau$ in this case and that this specific decomposition is the only one that uses no signals with $h_i^{-1}(1)\supseteq S^1$. (If $\bar f_*$ were larger, by Proposition~\ref{p-baseline-levels} there would be no such decomposition; if $\bar f_*$ were smaller, we will show later in Proposition~\ref{p-number-of-signals-pusc} that there would be multiple such decompositions.)

\begin{figure}[htbp]
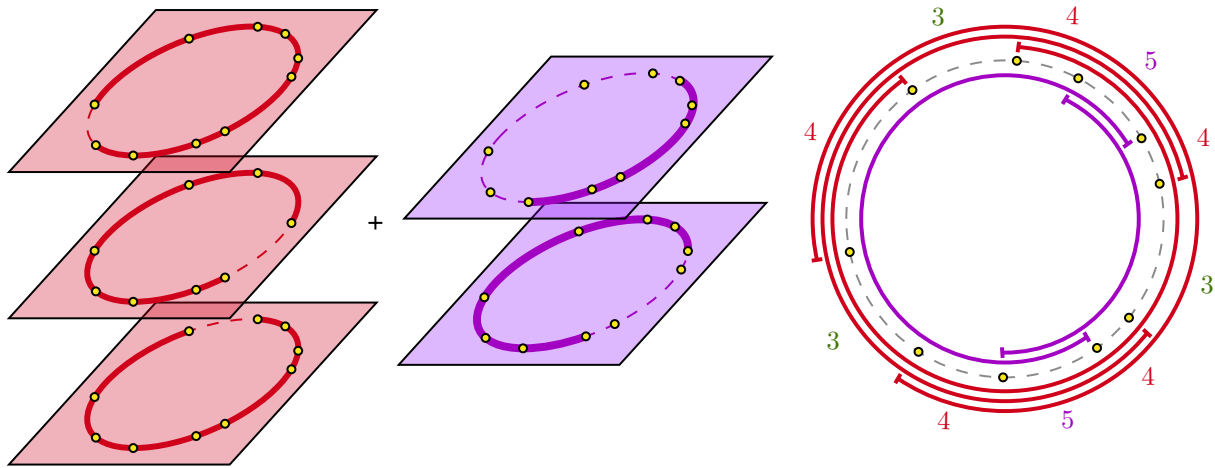


\centering

\tikzset{every picture/.style={line width=0.75pt}} 


\caption{A properly u.s.c.\ decomposition of the function from Example~\ref{x-three-sets} using no signals with $h_i^{-1}(1)\supseteq S^1$.}\label{f-no-baseline}
\end{figure}

\section{Degrees of freedom}\label{s-dof}

Now we want to show that all properly u.s.c.\ decompositions satisfying a suitable additional condition are equally likely in the sense that they maximize some measure of likelihood. For a decomposition $H_{n}(v)$ and an integer $N \ge n$, we define the $N$\emph{-degrees of freedom}
\begin{equation}\label{N-DOF}
    \text{DoF}_{N}(H_{n}(v)):=3(N-n)+\sum_{i=1}^{n}\rho( h_{i}),\qquad
\rho(h_{i}):= \begin{cases} 3 & \text{if } \card\bdry h_{i}=0   \\ 1 & \text{if } \card\bdry h_{i}\in\{1,2\}\\ 0 & \text{if } \card\bdry h_{i}=\infty\end{cases}\;.
\end{equation}
This definition alludes to a measure of likelihood because ``throwing" $N$ signal functions randomly into the plane, i.e., randomly choosing $(m_i,r_i)$, the dimension of the parameter space for $h_i$ with $\card\bdry h_i=0$ is 3, but when $\bdry h_i=\{x_1,x_2\}$ for some $x_1,x_2\in S^1$, this dimension gets restricted to 1 because we then know $m_i\in \R^2$ must lie in a line from the origin perpendicular to the line connecting $x_1$ and $x_2$, with an appropriate $r_i$. The case of $\bdry h_i=\{x_1\}$ is similar, only $m_i$ satisfying this condition lie on the line from the origin intersecting $x_1$, with an appropriate $r_i$. The case $\card\bdry h_{i}=\infty$ (i.e., $\bdry h_i=S^1$) corresponds to only one $h_i$, thus the dimension of the parameter space is 0. Furthermore, considering perturbations of $h_i$, we see that $\rho(h_i)$ is equivalent to the dimension of the subset of some sufficiently small neighborhood of $(m_i,r_i)$ such that $\bdry h_i$ stays exactly the same. For example, when $\card\bdry h_i=2$, we can perturb the radius and the center of the landmark in only one way (moving on the line described above), so that the boundary points of the perturbed decomposition stay the same. The term $3(N-n)$ in~(\ref{N-DOF}) is a correction term: we ``add" $N-n$ many $h_i$ with $\card\bdry h_i=0$ and $h_i^{-1}(1)\not\supseteq S^1$ to $H_n(v)$. These added $h_i$ are identically 0 on $S^1$, so the values of $H_n(v)$ are unchanged and the decomposition now has $N$ signals. The reason for doing this is to be able to use $\text{DoF}_N(\cdot)$ for comparing $H_{n_1}(v), H_{n_2}(v)\in\bar f$ with $n_1\neq n_2$. With the addition of this corrective term, we can deem the decomposition with a larger $\text{DoF}_N(\cdot)$ more likely. Additionally, we will see later that all robust decompositions of $\bar f $ will have the same $\text{DoF}_N(\cdot)$ under this formulation.

\begin{Proposition}\label{p-pusc-same-dof} All properly u.s.c.\ ${H_n(v)\in \bar f}$ with ${\card\bdry h_i\in \{0,2\}}\,\forall i$ have the same number of $N$-degrees of freedom, for every fixed $N\ge n$.
\end{Proposition}

\textbf{Proof.} By Corollary~\ref{c-intervals}, all properly u.s.c.\ $H_n(v)\in \bar f$ have the same number, say  $\ell$, of $h_i$ with $\card \bdry h_i=2$ in (\ref{coupledsignal}). By assumption, no $h_i$ with $\card\bdry h_i\in\{1,\infty\}$ exists in (\ref{coupledsignal}).  Any two properly u.s.c.\ decompositions $H_{n_1}(v_1),H_{n_2}(v_2)\in \bar f$ thus only differ in the counts of $h_i$ with $\rho(h_i)=3$; let these be $m_1$ and $m_2$, respectively.
For $i=1,2$, we have $n_i=\ell+m_i$ and so $\text{DoF}_N(H_{n_i}(v_i))=3(N-n_i)+3m_i+\ell=3(N-n_i+m_i)+\ell=3(N-\ell)+\ell$ which is the same for both decompositions.~\qed

\begin{Proposition}\label{p-pusc-max-dof} A decomposition $H_n(v)\in \bar f$ is properly u.s.c.\ and has $\card\bdry h_i\in \{0,2\}\,\forall i$ if and only if it maximizes $N$-degrees of freedom, for every fixed $N\ge n$.
\end{Proposition}

\textbf{Proof.} (if) We proceed by contraposition. Suppose first that $H_n(v)$ is not properly u.s.c. By Lemma~\ref{l-pusc-char}, $H_n(v)$ has either some $h_i$ with $\card\bdry h_i=1$ and $h_i^{-1}(1)\not\supseteq S^1$, or some $h_i,h_j$ whose supports on $S^1$ locally intersect only at a point. In the former case, such an $h_i$ has value 1 at only a single point in $S^1$. Therefore, the decomposition $H_{n-1}(v'):=H_n(v)- h_i$ would only disagree with $f$ on a finite number of points, i.e., it is still in $\bar f$. It is easy to check that $\text{DoF}_N(H_{n-1}(v'))=\text{DoF}_N(H_{n}(v))+2$, thus $H_n(v)$ does not maximize $\text{DoF}_N(\cdot)$. In the latter case, if the supports of $h_i,h_j$ only intersect at $x$ and possibly some other point, they could be replaced by any $h_i'$ with $\card\bdry h_i=2$ that has the set $(h_i^{-1}(1)\cup h_j^{-1}(1))\cap S^1$ as its support in $S^1$; otherwise if the supports of $h_i,h_j$ intersect in $S^1$ at some circular arc other that $x$, they could be replaced by some $h_i'$ with $\card \bdry h_i =2$ with support as this circular arc and some $h_j'$ with $\card\bdry h_j=0$ and $h_j^{-1}(1)\supseteq S^1$. This yields a valid decomposition of $f$ in both cases. Again such a replacement increases $\text{DoF}_N(\cdot)$ by 2 in both cases, therefore $H_n(v)$ does not maximize $\text{DoF}_N(\cdot)$.

Now suppose that $H_n(v)$ has some $\card\bdry h_i\in\{1,\infty\}$. We only need to consider $h_i^{-1}(1)\supseteq S^1$, because we are done by the previous argument otherwise.
But such an $h_i$ can be replaced by some $h_i$ with $\card\bdry h_i=0$ and  $h_i^{-1}(1)\supseteq S^1$, forming a valid decomposition of $f$ and increasing $\text{DoF}_N(\cdot)$ (by either 2 or 3). Therefore, $H_n(v)$ again does not maximize $\text{DoF}_N(\cdot)$.

\noindent (only if) By Proposition~\ref{p-pusc-same-dof}, given some $f$ we know that $N$-degrees of freedom is an invariant of all its properly u.s.c.\ decompositions with ${\card\bdry h_i\in \{0,2\}}\,\forall i$. By the ``if" part of the statement, it follows that all such decompositions maximize $N$-degrees of freedom.~\qed

From Propositions~\ref{p-pusc-max-dof} and~\ref{p-robust} we immediately obtain the following.

\begin{Corollary} Robust decompositions $H_n(v)\in \bar f$ maximize $N$-degrees of freedom, for every $N\ge n$.
\end{Corollary}

We can now also complete the discussion by providing a result justifying why the property of maximizing $N$-degrees of freedom is desired from a decomposition. We find that such decompositions have the least number of signals, in the sense of the following corollary.
\begin{Corollary}\label{c-maxdofminimal}
    Let $H_n(v)\in \bar f$. For every $H_m(v_m)\in \bar f$ that maximizes $N$-degrees of freedom, uses the same number of signal functions $h_i$ with $h_i^{-1}(1)\supseteq S^1$ as $H_n(v)$ does, and has no $h_i^{-1}(1)\cap S^1=\emptyset$, we have $m\leq n$.
\end{Corollary}

\textbf{Proof.} Take $H_n(v)$ and $H_m(v_m)$ in $\bar f$ as in the statement of the corollary. Since $H_m(v_m)$ maximizes $N$-degrees of freedom, by Proposition~\ref{p-pusc-max-dof} it is properly u.s.c.\ and only uses signals $h_i$ with $\card \bdry h_i\in\{0,2\}$. We know that $H_m(v_m)$ uses the
same number of signals $h_i$ with $h_i^{-1}(1)\supseteq S^1$ as $H_n(v)$ does, but it can differ from $H_{n}(v)$ in the number of signals $h_i$ with $\card\bdry h_i=2$ or $h_i^{-1}(1)\cap S^1=\emptyset$. By hypothesis, $H_m(v_m)$ does not use any signals of the latter type, and being properly u.s.c., by Corollary~\ref{c-intervals} and Remark~\ref{c-intervals_lowerbound} it uses the least possible number of signals of the former type among the decompositions in $\bar f$. In addition, $H_n(v)$ can have signals $h_i$ with  $\card \bdry h_i\in \{1,\infty\}$, but $H_m(v_m)$ does not. Thus $m\le n$, as claimed.~\qed

\section{Generating and counting robust decompositions}\label{s-generating-counting}

This section is concerned with generating and counting decompositions that maximize $\text{DoF}_N(\cdot)$ and, in particular, robust decompositions.
We start by noting that, if a decomposition $H_n(v)\in \bar f$ contains any signal functions that are identically 0 on $S^1$, i.e.,
$h_i$ with $\card \bdry h_i=0$ and $h_i^{-1}(1)\cap S^1=\emptyset$, then we can subtract them from $H_n(v)$ and the remaining decomposition will still be in $\bar f$. Accordingly, in what follows we ignore such ``inactive" signal functions when we do the counting.

\begin{Proposition}\label{p-number-of-signals-pusc} The number $n$ of signal functions $h_i$ in a properly u.s.c.\ decomposition ${H_n(v)\in \bar f}$ that maximizes $\text{DoF}_N(\cdot)$ satisfies
\begin{equation}\label{countbound}
    \max\{\bar f^*,\bar f^*-\bar f_*+\tau\}\leq n\leq \bar f^*+\tau \;.
\end{equation}
\end{Proposition}

\textbf{Proof.} Let $H_n(v)\in\bar f$ maximize $\text{DoF}_N(\cdot)$.
By Proposition~\ref{p-pusc-max-dof}, it has $\card \bdry h_i\in\{0,2\}$ $\forall i$. The number of $h_i$ with $\card \bdry h_i=2$ is fixed by Corollary~\ref{c-intervals} to be $\bar f^*-\bar f_*+\tau$.
It remains to count the $h_i$ with $\card\bdry h_i=0$, and among these, by the remarks immediately preceding the proposition, we are only interested in the ones with $h_i^{-1}(1)\supseteq S^1$ (because otherwise $h_i\equiv 0$ on $S^1$). Clearly, $H_n(v)$ cannot have more than $\bar f_*$ such $h_i$, which gives the upper bound in~\eqref{countbound}. On the other hand, we know from Proposition~\ref{p-baseline-levels} that $H_n(v)$ must have at least $\max\{\bar f_*-\tau,0\}$ such $h_i$. Thus we have $n\ge \bar f^*-\bar f_*+\tau+\max\{\bar f_*-\tau,0\}=\max\{\bar f^*,\bar f^*-\bar f_*+\tau\}$ which gives the lower bound in~\eqref{countbound}.~\qed

There is a close relationship between some of our findings---particularly the above Proposition~\ref{p-number-of-signals-pusc} and the earlier results in Section~\ref{s-top-prop}---and the results of~\cite{euler}. That paper also deals with counting signal functions based on the function given by their local sum. The setting of~\cite{euler} is more general than ours in several respects, one of which is that the functions are defined on a general topological space. To match our setting, we can take this space to be $S^1$. To be able to apply the results of~\cite{euler}, we need to assume that the supports of the individual signal functions, i.e., the sets $h_i^{-1}(1)\cap S^1$, all have the same nonzero Euler characteristic. In our case, this means that we can work with $h_i$ for which  $\card \bdry h_i=2$; then $h_i^{-1}(1)\cap S^1$ are contractible arcs and their Euler characteristic equals 1. By Theorem 3.2 of~\cite{euler}, the number of such $h_i$ is obtained by integrating $f$ over $S^1$ with respect to the Euler characteristic, and in our notation this number comes out to be $\bar f^*-\bar f_*+\tau$. This agrees with our Corollary~\ref{c-intervals}, and in fact our proof of Proposition~\ref{p-discont}---which eventually led to Corollary~\ref{c-intervals}---is based on a computation carried out in~\cite[Section~4]{euler} using so-called upper excursion sets (as we indicated in that proof). The other type of signal functions considered in Proposition~\ref{p-discont}, namely, the ones with $\card\bdry h_i=0$ and $h_i^{-1}(1)\supseteq S^1$, are not covered by the approach of~\cite{euler} because $S^1$ has Euler characteristic 0, and so we have to count them separately. We know, though, that even if such $h_i$ are present---and they must be present when $f_*>\tau$ by Proposition~\ref{p-baseline-levels}---the above count of $h_i$ with $\card \bdry h_i=2$ is still valid. Another difference with~\cite{euler} is that they assume the function $f$ to be given everywhere and do not work with equivalence classes $\bar f$ as we do. However, we know from Section~\ref{s-top-prop} that every equivalence class $\bar f$ contains a unique properly u.s.c.\ function, $\langle \bar f\rangle$, and the approach of~\cite{euler} can be applied to this function.

In summary, if we want to apply the results of~\cite{euler} to our problem of counting signal functions in possible decompositions $H_n(v)\in\bar f$ that maximize $\text{DoF}_N(\cdot)$, we can use the following procedure.\smallskip

\textbf{Procedure 1}
\begin{enumerate}
    \item Construct $\langle \bar f\rangle\in \bar f$ that is properly u.s.c. (which amounts to requiring~(\ref{properUSC}) at any discontinuity point).
    \item For $\max\{\bar f_*-\tau,0\}\le k\le \bar f_* $, let $f_k:=\langle \bar f\rangle-k$.
    \item Apply Theorem~3.2 of~\cite{euler} to each $f_k$, which gives $\bar{(f_k)}^*-\bar{(f_k)}_*+\tau(f_k)=\bar f^*-\bar f_*+\tau(\bar f)$ as the number of signals $h_i$ with $\card \bdry h_i=2$ (in decompositions using only such signals).
    \item For each $f_k$ with $k>0$ (if any), add $k$ to the number obtained at step~3 to arrive at the total number of $h_i$.
\end{enumerate}

The possible counts obtained using this algorithm agree with the bounds
of our Proposition~\ref{p-number-of-signals-pusc}.

We now want to show the existence of $H_n(v)\in \bar f$ for each $n$ satisfying the bounds~\eqref{countbound} of the Proposition~\ref{p-number-of-signals-pusc}.
To illustrate, let us consider again $\bar f$ from Example~\ref{x-three-sets}. Applying Proposition~\ref{p-number-of-signals-pusc}, we see that the number of signal functions must satisfy $5\leq n \leq 8$. A decomposition with $n=8$ was already shown in Figure~\ref{f-example-3}, and one with $n=5$ was shown in Figure~\ref{f-no-baseline}. We now complete the picture by providing (non-unique) decompositions with $n=6$ and $n=7$ in Figure~\ref{f-n=6-7} below.
It is useful to observe that they utilize the same sets as the decompositions shown in Figures~\ref{f-example-3} and~\ref{f-no-baseline}, but combine them in different ways.

\begin{figure}[htbp]
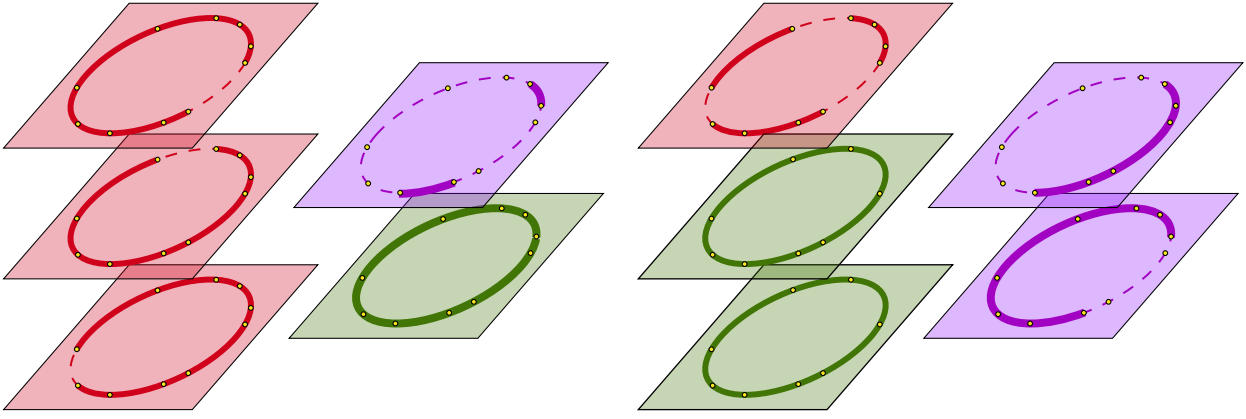

\centering


\caption{\small Two more properly u.s.c.\ decompositions of $H_n(v)$ of the function from Example~\ref{x-three-sets}: one with $n=6$ including one signal with $h_i^{-1}(1)\supseteq S^1$ (left), and one with $n=7$ including two signals with $h_i^{-1}(1)\supseteq S^1$ (right).}
\label{f-n=6-7}
\end{figure}

Looking at these examples, Figure~\ref{f-no-baseline} suggests a way of generating minimal decompositions for any $\bar f$---in the sense that the decomposition uses the least number (given by Proposition~\ref{p-baseline-levels}) of $h_i$ with $h_i^{-1}(1)\supseteq S^1$. For now, let $\bar f_*\geq \tau$. To generate such a decomposition of $\bar f$, we can use the following procedure.\smallskip

\textbf{Procedure 2}
\begin{enumerate}
    \item Consider the chain of upper excursion sets~\eqref{upperexcursionchain} for $\langle \bar f\rangle$ (Figure~\ref{f-example-3} might help visualize).
    \item For each $\bar f_*+1 \leq k\leq \bar f ^*$: For each connected component $C_i$ of $S^1\setminus \langle \bar f  \rangle^{-1}([k,\infty)) $, construct some signal function $h_i$ with $\card \bdry h_i=2$ and has support $S^1\setminus C_i$ in $S^1$. Summing each $h_i$, we generate a decomposition $H_{\tau_k+1}(v_k)$ with $\min _{S^1}(H_{\tau_k+1}(v_k))=\tau_k$ and $H_{\tau_k+1}(v_k)^{-1}([\tau_k+1,\infty))=\langle \bar f  \rangle^{-1}([k,\infty))$.


    \item Add $\bar f_*-\tau$ many signal functions $h_i$ with $h_i^{-1}(1)\supseteq S^1$ to $\sum _{k=\bar f_*+1}^{\bar f^*}H_{\tau_k+1}(v_k)$ to obtain a properly u.s.c. decomposition in $\bar f$.
\end{enumerate}

\begin{Proposition}\label{p-existence}
     $H_n(v)\in \bar f$ that maximizes $N$-degrees of freedom exists for all $n\in \mathbb{N}$ satisfying the bounds~\eqref{countbound}. $\omitted$
\end{Proposition}

Proposition~\ref{p-pusc-char} states that for properly u.s.c.~decomposition, each $x\in S^1$ belongs to $k_x$ many signal functions $h_i$ with $\card \bdry h_i=2$. Furthermore, the equivalent condition on the statement of Lemma~\ref{l-pusc-char} requires that such $h_i$ must have support in the same direction around $x$, i.e., it requires that they don't locally intersect only on a point around $x$.  Hence for any discontinuity point of properly u.s.c. $H_n(v)\in\bar f$, the support---in the shape of an arc of the circle---of the signal function $h_i$ with $\card\bdry h_i=2$ always lies in the direction in which $\langle\bar f \rangle$ increases.

The class of decompositions $H_n(v)$ we want to consider maximizes $N$-degrees of freedom---is properly u.s.c.~and uses robust signal functions $h_i\, \forall i$. Given some $\bar f$ generated from agent data, we want to find a complete set of all such ${H_n(v)\in\bar f}$. We can obtain a neat way of counting using the above fact that each discontinuity point of $\langle\bar f\rangle$ dictates a count and direction for the support of signal function(s) $h_i$ with $\card\bdry h_i=2$. The proof of the following Proposition is just the application of this idea with the necessary bookkeeping.

\begin{Proposition}\label{p-counting}
If some robust $H_n(v)\in \bar f$ exists and $\bar f_* \geq \tau$, then there are $(\bar f^*-\bar f_*+\tau)!$
unique robust ${H_n(v)\in\bar f}$. $\omitted$
\end{Proposition}

We revisit $\bar f$ from Example~\ref{x-three-sets}. $\bar f _*=\tau$, and a robust decomposition exists---e.g. Figure~\ref{f-example-3},~\ref{f-n=6-7},~\ref{f-no-baseline}. Hence Proposition~\ref{p-counting} applies, and there are $(\bar f ^*-\bar f_*+\tau)=(5-3+3)!=120$ unique decompositions.

In general, when $\bar f_*\geq \tau$, the problem of counting all properly u.s.c.\:$H_n(v)\in \bar f$ that maximize $DoF(\cdot)$ by way of connecting directed boundary points is equivalent to counting all $m \times n$ \textit{bipartite graphs} subject to the constraint $k_x$, dictating the number of edges to be connected to the node $x\in \bdry H_n(v)$. Interested reader might consult~\cite{bipartite} for specific cases.

When $\bar f_*<\tau$, this strategy produces invalid decompositions with just the sum of signal functions $h_i$ with $\card \bdry h_i=2$ having a minimum larger that $f_*$ on $S^1$. Accordingly, some decompositions must be rejected, and deciding is not easy. In this case, similarly as above, we can only give bounds for the number of decompositions.
\medskip

\noindent\textbf{Acknowledgment.} The authors are indebted to Yuliy Baryshnikov and Sayan Mitra for stimulating technical discussions.

\bibliography{robustdecompositions}
\bibliographystyle{siam}
\end{document}